\documentclass[11pt]{article}
\usepackage{amsmath}
\usepackage{setspace}
\usepackage{fullpage}
\usepackage{graphicx}
\usepackage{amsthm}
\newtheorem{theorem}{Theorem}
\newtheorem{lemma}{Lemma}
\newtheorem{definition}{Definition}
\newtheorem{proposition}{Proposition}
\usepackage{amsmath}            %important package for math     
\usepackage{epsf}               %package for including figures
\usepackage{amsfonts}
\begin{document}

\title{On The Number of Genus One Labeled Circle Trees}

\author{
Karola M\'esz\'aros\\
Massachusetts  Institute of Technology\\
 {\tt karola@math.mit.edu}
\\
%\vspace{1in}
}

\date{}

\maketitle

\begin{abstract}

\begin{small}

 A  genus one labeled circle tree is a tree with its vertices on a circle, such that
together they can be embedded in a surface of genus one, but not of genus zero.
We define an e-reduction process whereby a special type of subtree, called an
e-graph, is collapsed to an edge. We show that genus is invariant under
e-reduction. Our main result is a classification of 
genus one labeled circle trees through e-reduction. Using this we prove  a modified version of a conjecture of David Hough, namely,  that the number of
genus one labeled circle trees on $n$ vertices is divisible by $n$ or if it is not divisible by $n$ then it is divisible by $n/2$. Moreover, we explicitly characterize when each of these possibilities occur.

%\textbf{Keywords:} labeled circle tree, genus one circle tree 
\end{small}
\end{abstract}

\begin{singlespace}

%\vspace{-0.1in}

\section{Introduction}

Graphical enumeration arises in a variety of contexts in combinatorics \cite{ge}, and naturally so in the realm of combinatorial objects with interesting topological properties  \cite{gs}.  We provide a new classification of genus one circle trees and  address a question raised by Hough  \cite{hou} about their number. Our study is motivated by numerous results in the study of partitions and trees of a certain genus, as well as results about the genuses of  maps and hypermaps,   \cite{coma}, \cite{noy},  \cite{si}, \cite{wl}.

The following two definitions are discussed in \cite{hou} in great detail;  we shall use the definition of a labeled circle tree throughout the paper, whereas we shall mostly use an alternate, less technical definition for the genus of a circle tree. 

\begin{definition}
 A \textbf{labeled circle tree} (\textbf{l-c-tree}) on $n$ points is a tree with its 
  $n$ vertices labeled  1 through n on a circle in a
  counterclockwise direction
  and its edges  drawn as straight lines within the circle. 
\end{definition}

\begin{definition}  
The \textbf{genus} g \textbf{of a l-c-tree $T$} on $n$ points
is defined to be $g(\alpha)=$1+$\frac{\it{1}}{\it{2}}$(n-1-$\it{z}(\alpha)-\it{z}({\alpha}^{-1} \cdot \sigma))$,
where $\alpha$ is the matching  of the given l-c-tree $T$,
 $\sigma=$(1 2 3 4 $\ldots$ n), and the function $\it{z}$ gives the number of cycles
 of its argument; $\alpha^{-1}$ is the inverse permutation of $\alpha$,
 and the multiplication of two permutations is from right to left.
 \end{definition}

The genus of  a l-c-tree can also be described as the genus of the surface with
minimal genus such that the tree together with the circle it is
drawn on can be drawn on the surface without crossing edges. In particular, 
a genus one l-c-tree is such that the tree together with the circle it is drawn on can be embedded in a surface of genus one, but not of genus zero.

 Hough  \cite{hou} observed that 
the number of genus one labeled circle trees on $n$ points (denoted by \textbf{\textit{f(n)}}) is divisible by $n$ for small values of $n$, and hypothesized the same for 
  all integers $n>3$.  Using our classification of all genus one labeled circle trees,  we  prove that either $f(n)$ is divisible by $n$, or it is divisible by $\frac{n}{2}$; moreover, we explicitly describe when each of these possibilities occur.

In Section 2 we discuss the necessary definitions and review a result of  Marcus \cite{mar}, which implies that deleting an uncrossed edge from a l-c-tree or deleting all but one of several parallel edges leads 
  to one of two canonical reduced forms of circle trees 
 if and only if the l-c-tree was genus one. 
Although Marcus' result \cite{hou} is formulated for partitions, it easily translates to l-c-trees. Since the labeling of the l-c-tree is irrelevant for the deletion of edges mentioned above, we introduce the concept of an \textbf{\textit{unlabeled circle tree}} to which Marcus' result still applies.  For understanding the interrelation between the number of genus one l-c-trees and genus one u-c-trees  we explore the basic properties of u-c-trees in Section 3. We  introduce a special-structured subgraph, called an \textbf{\textit{edgelike-graph}}, in Section 4,  and  we describe an \textbf{\textit{e-reduction  process}} in Section 5. Based on the e-reduction process  we give a \textbf{\textit{classification}} of genus one c-trees by nineteen \textbf{\textit{reduced forms}} in Section 6.  We  clarify the connection between the number of l-c-trees and u-c-trees on $n$ points in the further sections, analyzing  reduced forms. Finally, using our previous results   we  formulate the theorem about $f(n)$'s divisibility by $n$ or $\frac{n}{2}$.

\section{Definitions and Remarks}

The definition of a labeled circle tree straighforwardly extends to the definition of a \textbf{\textit{labeled circle graph}}. Indeed, replacing the word tree with graph in the definition of l-c-tree gives the desired definition of a l-c-graph. Two \textbf{\textit{l-c-graphs}} $G_1$ and $G_2$ are said to be \textbf{\textit{isomorphic}}   
if an edge $e_1$ with endpoints labeled $i$ and $j$ is in $G_1$ if and only if there is an 
edge $e_2$  in $G_2$ with endpoints $i$ and $j$ (we consider graphs without multiple edges). Furthermore, if a vertex labeled $k$ is of degree zero in one of the graphs then it is of degree zero in both of them. 
An \textbf{\textit{unlableled circle graph}} (\textbf{\textit{u-c-graph}}) is a graph obtained by deletion of labels of a l-c-graph. Two \textbf{\textit{u-c-graphs}} are said to be \textbf{\textit{isomorphic}} if it is possible to label their vertices so that the obtained l-c-graphs are isomorphic.   

It follows from the definition of genus that isomorphic l-c-trees have equal genuses. We conveniently define the \textbf{\textit{genus of a u-c-tree}} $T$ on $n$ points  to be the genus of any of the l-c-trees one obtains by labeling the vertices of $T$ by \textit{1} through $n$ in a counterclockwise direction.

% Note that the genus of a u-c-tree is well defined since the genuses of all the possible l-c-trees obtained by the process described are equal. 

Call  u-c-graphs (u-c-trees) and l-c-graphs (l-c-trees)  by the common name \textbf{\textit{c-graphs}} (\textbf{\textit{c-trees}}).   
 \textbf{\textit{Edges}} $e_1$ and $e_2$ of a c-graph  \textbf{\textit{cross}} if
they have a point in
common on the drawing of the c-graph  other than their endpoints.
\textbf{\textit{Edges}} $e_1$ and $e_2$ of a c-graph $G$  are \textbf{\textit{parallel}} if 
they cross the same edges of $G$, respectively.
That the  relation `parallel' is an equivalence relation is a straightforward check of reflexivity, symmetry and transitivity.

\begin{definition}
 An \textbf{u-c-tree} C and a \textbf{l-c-tree} T are said to  \textbf{correspond}  if the u-c-tree obtained by the deletion of the labels of T  is isomorphic to C. 
\end{definition} 

By the definition of genus the genuses of corresponding u-c-trees and l-c-trees are equal. 
%\begin{center}
%\epsfbox{tree.eps}
% \end{center}
Note that   
 a u-c-tree $C$ on  $n$ points can correspond to at most $n$  non-isomorphic l-c-trees. In some cases
 the u-c-tree corresponds to exactly $n$ non-isomorphic l-c-trees, 
  but in some cases a u-c-tree corresponds to less then $n$ non-isomorphic 
   l-c-trees,  Figure 2.1. 
 %\begin{center}

%\epsfbox{ctree11.eps} 
 
%\epsfbox{ctree12.eps} 
\begin{figure}
\begin{center}
\epsfbox{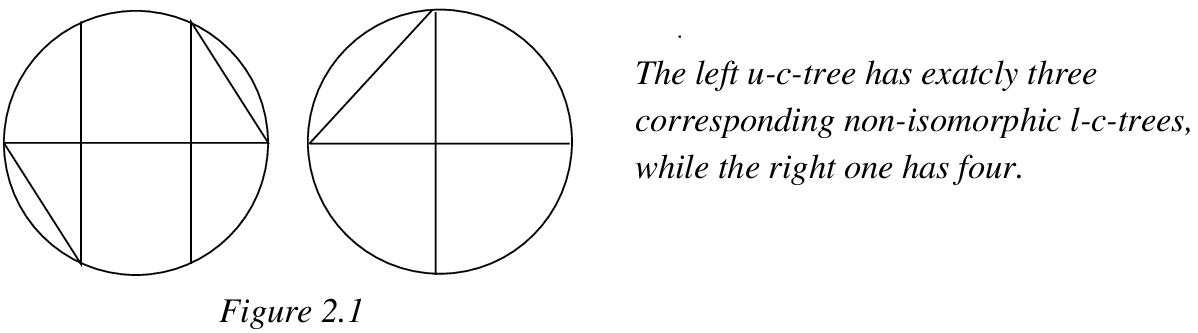}
\end{center}
\end{figure}

The reinterpretation of Marcus' result [3, 4]:

\begin{proposition} Performing the following two operations on a c-tree as many times as possible: 
 
1) deleting an edge from the c-tree, which is not crossed by any other 
edge 
 
2) deleting all but one of several parallel edges

 leads to Form 1 or Form 2  
shown on Figure 2.2 if and only if the c-tree was genus one.
\label{p1} 
\end{proposition}

\begin{figure}
 \begin{center}
\epsfbox{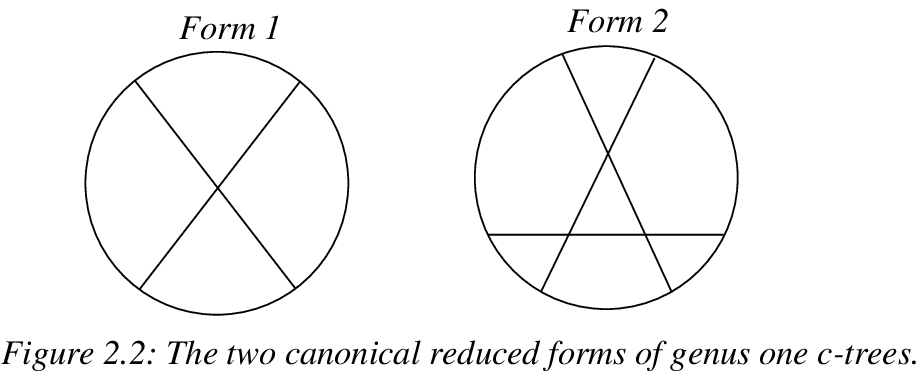} 
 \end{center}
\end{figure}

We refer to the two operations of Proposition 1 as operation 1) and operation 2).

\begin{definition} Call the u-c-graphs obtained from a u-c-tree C
by executing operations 1) and 2)
\textbf{offsprings}. A u-c-tree $C$ \textbf{descends} to a u-c-graph, if the u-c-graph is an offspring of $C$.
The \textbf{final offspring} of a u-c-tree C is the offspring which has no edges which could be deleted by the execution of operations 1) and 2). (Some offsprings of a u-c-tree are represented on Figure 2.3)
\end{definition}

\begin{figure}
\begin{center}
\epsfbox{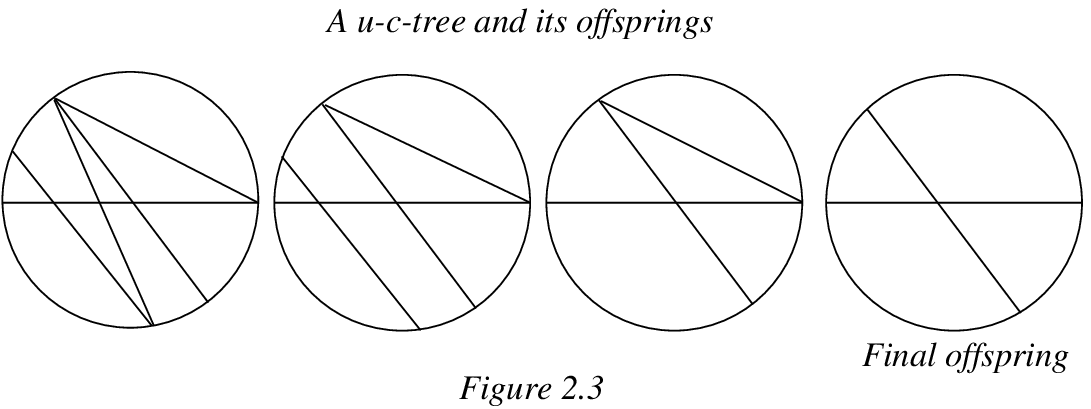}
\end{center}
\end{figure}

In Section 3 we conclude  that the final offspring of a genus one  u-c-tree C is unique (up to isomorphism).

We now rephrase Proposition 1:

$\bf{Proposition}$ $\bf{ 1'.}$ 
\textit{The genus of a u-c-tree T is one if and only if its final offspring is Form 1 or Form 2.} 

Naturally,  by saying that the final offspring is Form 1 or Form 2, we  mean that the final offspring is isomorphic either to Form 1 or to  Form 2. We do not  stress this in the future, since it is clear from the context.
 
We can now reformulate  the question about the divisibility of  $f(n)$ by $n$ or $\frac{n}{2}$ as follows: When is the number of l-c-trees on $n$ points, which have corresponding   u-c-trees that
 descend to Form 1 or Form 2  (Figure 2.2) by the execution of operations 1) and 2) (these are all of the genus one l-c-trees), divisible by $n$ and when is it divisible just by $\frac{n}{2}$ and not by $n$?

\section{Initial Observations Concerning U-C-Trees}

We state two simple lemmas concerning u-c-trees without proof. The proofs are based on the definitions of uncrossed and parallel edges, and are easily derived by contracition. 
  
\begin{lemma} 
 
If an edge is uncrossed after 
a number of operations 1) and 2) are executed on a u-c-tree, then that edge 
is uncrossed in the u-c-tree, as well as in all its offsprings (if not deleted). 
\end{lemma} 

\begin{lemma} 
If two edges are parallel after 
a number of operations 1) and 2) are executed on a u-c-tree, then those two edges 
are parallel in the u-c-tree, as well as in all its offsprings (if not deleted). 
 \end{lemma} 

 From Lemma 1 and Lemma 2 we conclude
 that the order of the execution of operations 1) and 2) and the particular
 choice of the order of the edges to be deleted do not affect the final offspring. 
  By Proposition 1, after 
        executing operations 1) and 2) on a u-c-tree until applicable, 
        Form 1 or Form 2 are obtained if and only if the u-c-tree was genus one.
        Thus, it is possible to 
        construct every genus one u-c-tree by beginning from  Form 1 or Form 2 
 and  by adding parallel edges to the ones presented in the form, 
        and by adding uncrossed edges.  Moreover, starting with these two forms and adding only parallel and uncrossed edges any u-c-tree obtained is   of genus one. 
See Figure 3.1 for illustration. 
This ``building'' idea might  serve as a basis for obtaining the exact number of genus one l-c-trees on $n$ points.
%of the ``building process'' developed in Section 6.  

\epsfbox{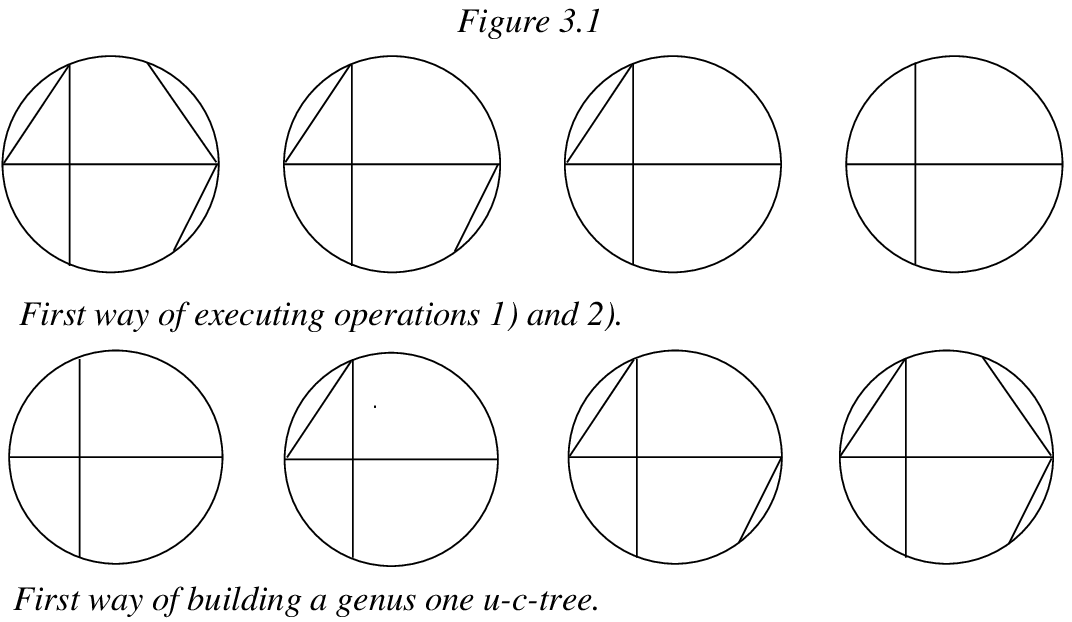}  

\epsfbox{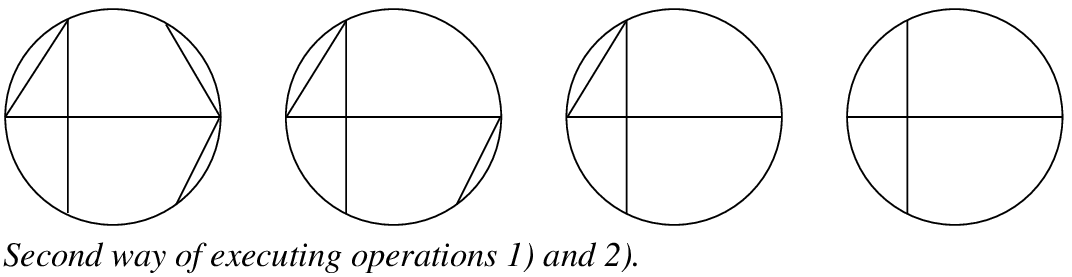} 

\epsfbox{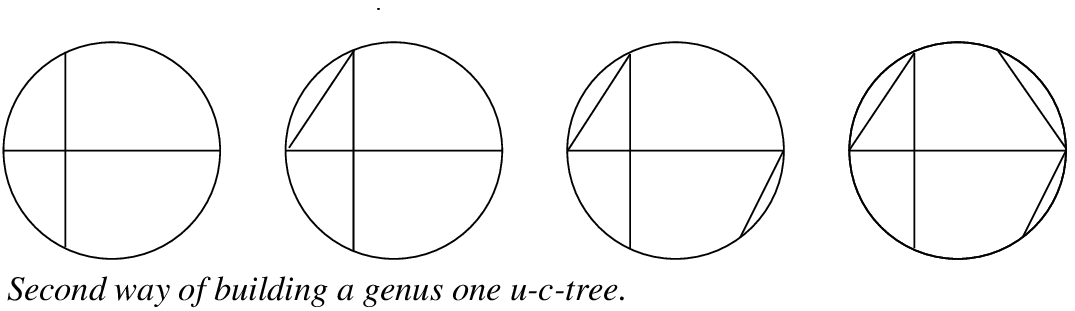}

\section{About Edgelike-Graphs}

In this section we introduce a main concept of our work, that of an \textbf{\textit{edgelike-graph}}. As the name already suggests, these graphs behave somewhat like edges. Indeed, edgelike-graphs, or e-graphs for short, are subtrees of a given u-c-tree with the special property that collapsing an e-graph to an edge (a specified one) the obtained u-c-tree has the same genus as the one we began with. 

 Once the concept of e-graph is grasped, the way operations 1) and 2) act on a c-tree becomes easy to visualize and understand.  A u-c-tree can be decomposed  into e-graphs, in which case  operations 1) and 2) act within these decomposed structures. The previous fact exhibits the correlation between the structure of a c-graph and operations 1) and 2). 

Let edges $e_1, e_2, \cdots, e_k$ be parallel.  If  $e_i$ and $e_j$ are the outermost edges among  $e_i, e_{i+1}, \cdots, e_{j}$,  for all $1\leq i\leq j\leq k$, then edges $e_1, e_2, \cdots, e_k$ are \textbf{\textit{increasingly parallel}}. Edges $CD$, $CG$, $JH$, $FE$ are increasingly parallel on Figure 4.1.
Edges
$e_1, e_2$, $\ldots$, $e_k$ constitute a \textbf{\textit{path}} if and only if
 there exist points
$E_1, E_2,\ldots, E_{k+1}$  on the circle such that the endpoints of
$e_i$ are $E_i$ and $E_{i+1}$ for all $1 \le i \le k$.
We also make the convention that arc $\widehat{AB}$ is the arc between points $A$ and $B$ when going  in a counterclockwise direction from $A$ to $B$.      

\begin{definition}
Given a c-graph $G$ take any crossed edge $AB$ of it. Let increasingly parallel edges $e_1, e_2, \cdots, e_k$ be all edges of $G$ parallel to $AB$ (including $AB$ itself).  
Let increasingly parallel edges $a_1, a_2, \cdots, a_l$ $\in \{e_1, e_2, \cdots, e_k \}$ be all of the edges parallel to $AB$ such that there exist a path of edges $AB=b_1$, $b_2, \cdots, b_m=a_i$ ($i \in [l]$) such that each $b_j$, $j \in [m]$, is either uncrossed or parallel to $AB$. If $a_1=CD$ and $a_l=EF$, then the edges of $G$ such that both of their endpoints are on arcs $\widehat{DE}$ and $\widehat{FC}$ and they are uncrossed or parallel to $AB$ constitute the edgelike-graph, or e-graph, of $G$ containing $AB$. Figure 4.1.      
\end{definition}

Observe that there is a unique e-graph containing each crossed edge.

\begin{definition}
Let the arcs  $\widehat{DE}$ and $\widehat{FC}$ as in the above definition be called the arcs of an e-graph, whereas edges $CD$ and $EF$ the outermost edges of it. Also, call the set of crossed edges of the e-graph the set of parallel edges of the e-graph and call the set of edges of the e-graph that are uncrossed the  set of uncrossed edges. 
\end{definition}  

 \begin{figure} 
\begin{center}
\epsfbox{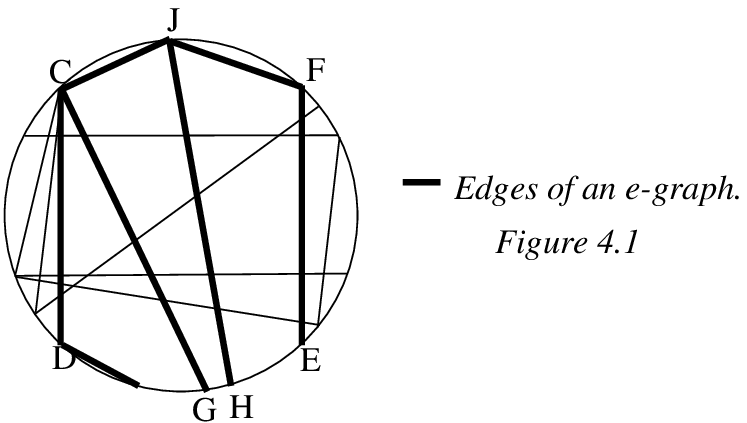}

%\epsfbox{el.eps} 
\end{center}
\end{figure}

\begin{lemma}
Using the notation of Definition 5, edges  $a_1, a_2, \cdots, a_l$ are all of the edges parallel to $AB$ having both of their endpoints on arcs $\widehat{DE}$ and  $\widehat{FC}$.
\end{lemma} 

\begin{proof}
Note that if $AB=e_i$, and $e_j=a_z$ for some $z\in [l]$, then any $e_r$, $r$ between $i$ and $j$, is equal to some $a_q$ for some $q \in [l]$. This observation leads to the proof of the lemma.   
\end{proof}

\begin{lemma} If $\mathcal{E}$ is an e-graph of c-graph $G$, then 
 $\mathcal{E}$ consists of edges having no points in common except for their vertices.
 \end{lemma}  

\begin{proof}
Suppose the opposite. Let $e_1$ and $e_2$ be two edges   of  $\mathcal{E}$
  having a point $A$ in common, such that $A$ is not
their endpoint.
Since any edge of  $\mathcal{E}$ is either uncrossed
 or parallel to an edge $e$ (being uncrossed and parallel considered within $G$), we 
conclude that $e_1$ and $e_2$ are both parallel to $e$. However, crossing edges cannot be parallel. 
  This contradiction proves the statement.
  \end{proof}

%\begin{lemma}

\begin{lemma}
If $\mathcal{E}$ is an e-graph of a u-c-tree C with 
arcs $\widehat{DE}$ and $\widehat{FC}$ and outermost edges $CD$ and $EF$, then there is no edge of C having one of its endpoints of the open arcs $\widehat{DE}$ or $\widehat{FC}$ and the other endpoint on the open arcs  $\widehat{CD}$ or $\widehat{EF}$. 
\end{lemma}

\begin{proof}
The statement of the lemma follows since $CD$ and $EF$ are parallel edges.
\end{proof}

\begin{proposition}
Given an e-graph $\mathcal{E}$ of a genus one u-c-tree C, let $KL$ and $MN$ be its outermost edges, and the arcs  $\widehat{LM}$ and $\widehat{NK}$ its arcs. Then all edges of C  which have both of their endpoints
on arcs $\widehat{LM}$ and $\widehat{NK}$ are edges of $\mathcal{E}$. Conversely, only
such edges are edges of  an e-graph $\mathcal{E}$ of a genus one u-c-tree C.
\end{proposition}

\begin{proof}
Let $\mathcal{U}$ and $\mathcal{P}$ be the sets of uncrossed and parallel edges of e-graph $\mathcal{E}$ 
 described in the proposition.
By the definition of e-graph $\mathcal{U}$ contains all the uncrossed edges of C with endpoints
on arcs $\widehat{LM}$ and $\widehat{NK}$ and $\mathcal{P}$ contains all edges of C  parallel to $KL$ with endpoints
on arcs $\widehat{LM}$ and $\widehat{NK}$. To prove  Proposition 2
it suffices to show that 
there is no edge $e$ of C with both of its endpoints on arcs $\widehat{LM}$ and $\widehat{NK}$
 such that it is not in $\mathcal{U}$ or $\mathcal{P}$. 
 Suppose the opposite, that there was an edge $e$ of C with both of its endpoints on arcs $\widehat{LM}$ and $\widehat{NK}$
 such that it is not in $\mathcal{U}$ or $\mathcal{P}$.  
If $e$ had one of its endpoints
 on $\widehat{LM}$ and the other on $\widehat{NK}$, then all the edges crossing KL would cross $e$.
 Since $e$ was not parallel to $KL$ there would have been some edge $e'$ which
  does not cross $KL$ and $MN$ but crosses $e$. 
 Edge
 $e'$ could clearly not be in $\mathcal{U}$, and  it also could not be in $\mathcal{P}$, since  $KL$ and
 $MN$ could not be crossed by $e$,  given that the
 endpoints of $e$ are on arcs $\widehat{LM}$ and $\widehat{NK}$. 
Edge $e'$ could be an edge
  with both endpoints on one of the arcs $\widehat{LM}$ or $\widehat{NK}$ or with one endpoint on $\widehat{LM}$
   and other endpoint on $\widehat{NK}$ (by Lemma 5 these are the only possibilities),  Figure 4.2. 
It is clear that executing
   operations 1) and 2) we would not
   get to Form 1 or Form 2, since the cross from $e$ and $e'$ and
   from $KL$ and some edge it crosses would remain.
\begin{figure}
\begin{center}
\epsfbox{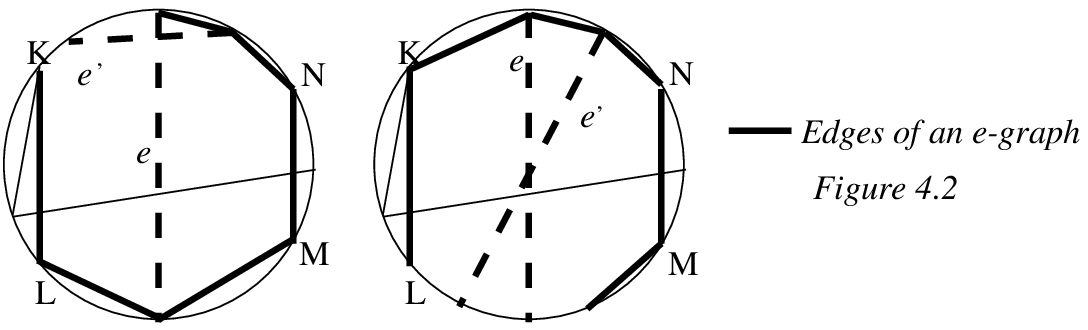}
\end{center} 
\end{figure}

On the other hand, if $e$ had both of its endpoints on one of the arcs $\widehat{LM}$ or $\widehat{NK}$, since it was not uncrossed it would have been crossed by some edge $e'$ and by Lemma 5  $e'$ would have both of its endpoints on arcs  $\widehat{LM}$ and $\widehat{NK}$. All cases are depicted on Figure 4.3. The cross obtained from the crossing of $e$ and $e'$ and the cross from $KL$ and some edge it crossed it would necessarily remain after executing operations 1) and 2) so we could not  get to Form 1 or Form 2, thus the genus of the u-c-tree could not be one.     

\begin{center}
\epsfbox{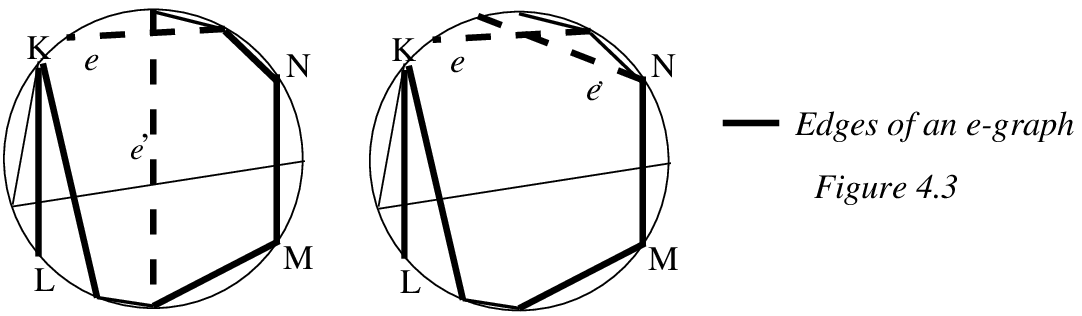}
\end{center} 

Thus,  all edges of C  which have both of their endpoints
on arcs $\widehat{LM}$ and $\widehat{NK}$ are edges of the e-graph  $\mathcal{E}$. It follows by definition that only such edges are edges of the e-graph. 

\end{proof}

\textbf{Corollary.} (Definition 5, Lemma 5, Proposition 2) \textit{An e-graph of a genus one u-c-tree is a tree.} 
Proposition 2 also implies that given the outermost edges of an e-graph, the e-graph is uniquely determined.

\begin{lemma}
Given two  e-graphs $\mathcal{E}_1$ and $\mathcal{E}_2$ in a genus one u-c-tree C
with sets of parallel edges $\mathcal{P}_1$ and $\mathcal{P}_2$, if $\mathcal{P}_1=\mathcal{P}_2$, then 
$\mathcal{E}_1$ and $\mathcal{E}_2$ are not different. (Two
e-graphs of a c-tree are said to be different if there is an edge in one of them
which does not belong to the other. When two e-graphs are not different, we also say that they are identical.)
\end{lemma}
\begin{proof}
The set of parallel edges of an e-graph determine its outermost edges and the outermost edges determine the e-graph in a genus one u-c-tree. 
\end{proof}

\begin{proposition}
There are no two different e-graphs $\mathcal{E}_1$ and $\mathcal{E}_2$ 
in a genus one u-c-tree C  such that they have vertices in common.
\end{proposition}

\begin{proof}
Let $\mathcal{U}_1, \mathcal{P}_1, \mathcal{U}_2, \mathcal{P}_2$ be the sets of uncrossed and parallel edges of
two different e-graphs $\mathcal{E}_1$ and $\mathcal{E}_2$. From Lemma 6 
$\mathcal{P}_1\ne \mathcal{P}_2$. If the edges of $\mathcal{P}_1$ and $\mathcal{P}_2$
are parallel it is impossible that the e-graphs have common vertices by the definition of an e-graph. 
In case the edges of $ \mathcal{P}_1$ and $ \mathcal{P}_2$ are not parallel, the only vertex two e-graphs
$\mathcal{E}_1$ and $\mathcal{E}_2$ might have in common is an endpoint
of some of their outermost edges. However, if 
$\mathcal{E}_1$ and $\mathcal{E}_2$ had such a point in common, there must have been some edge $e$ which crosses, say the edges of $ \mathcal{P}_2$ and does not  cross the edges of $ \mathcal{P}_1$. In this case the cross made by $e$ and some edge of $ \mathcal{P}_2$ as well as some cross from some edge of  $ \mathcal{P}_1$ and another edge, not parallel to $e$, must stay after executing operations 1) and 2) thus it is impossible to obtain Form 1 or Form 2 (Figure 4.4).
\begin{figure}
\begin{center}
\epsfbox{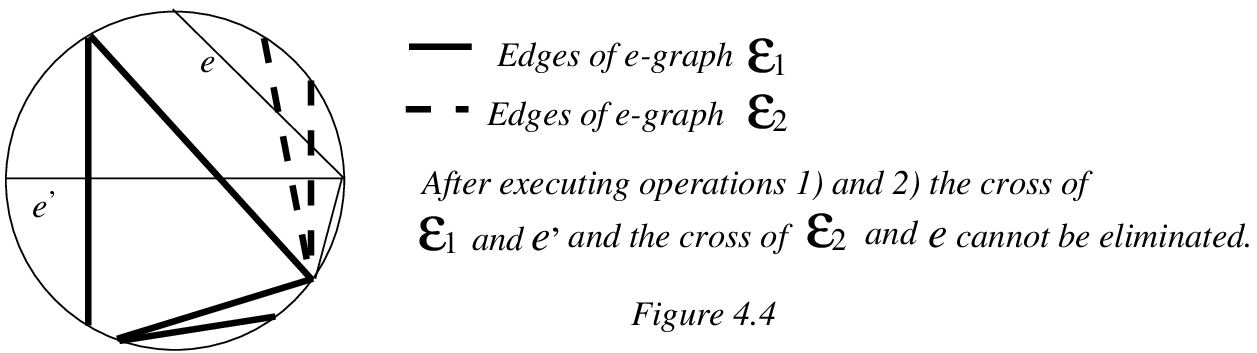}
\end{center}
\end{figure}
Thus, if the u-c-tree is genus one then no two e-graphs of
the u-c-tree have vertices in common.
\end{proof}

Given an  e-graph $\mathcal{E}$ of a u-c-graph with  its set of uncrossed edges $\mathcal{U}$
and  set of  parallel edges $\mathcal{P}$, call the elements
 of $\mathcal{P}$ the \textbf{\textit{parallel edges}} and the elements
 of $\mathcal{U}$ the \textbf{\textit{uncrossed edges}}
 of the e-graph $\mathcal{E}$.
We say that an e-graph $\mathcal{E}$ is \textbf{\textit{parallel}} to an edge $e$ if its parallel edges are parallel to $e$. Similarly,   e-graph $\mathcal{E}_1$ is parallel to another e-graph $\mathcal{E}_2$ if their parallel edges are parallel.  
We  say that an edge $e$ is \textbf{\textit{between}}  parallel e-graphs $\mathcal{E}_1$ and  $\mathcal{E}_2$  if both endpoint of $e$ are on arcs $\widehat{BE}$ and $\widehat{HC}$, where edges $AD$, $BC$, $EH$, $FG$ are increasingly parallel and the arcs of $\mathcal{E}_1$ are $\widehat{AB}$, $\widehat{CD}$ and the arcs of $\mathcal{E}_2$ are $\widehat{EF}$, $\widehat{GH}$ (and if $e$ is not an edge of $\mathcal{E}_1$ or  $\mathcal{E}_2$). For example, taking e-graphs $\mathcal{E}_1$  and $\mathcal{E}_2$ from Figure 4.5, an edge is between these two e-graphs if and only if both of its endpoints are on arcs $\widehat{EF}$ and $\widehat{GH}$.
A \textbf{\textit{path}} consisting of edges  $E_1E_2$, $E_2E_3$, $\ldots$, $E_kE_{k+1}$,  \textbf{\textit{connects}} 
two \textbf{\textit{e-graphs}} $\mathcal{E}_1$ and $\mathcal{E}_2$ if and only if point 
$E_1$ is the intersection of $\{E_1, E_2$,$\ldots$,$ E_{k+1}\}$ and the
points of one of the e-graphs, and point
  $E_{k+1}$ is the intersection of $\{E_1, E_2,\ldots, E_{k+1}\}$
  and the points of the other of  the e-graphs.

\begin{theorem} 
There can be at most two different e-graphs $\mathcal{E}_1$ and $\mathcal{E}_2$
of a genus one u-c-tree  C   such that $\mathcal{E}_1$ and    
$\mathcal{E}_2$ are parallel.     
\end{theorem}

\begin{proof}    
The main idea of the proof is that a u-c-tree is connected, and if there were already three parallel e-graphs the u-c-tree, then it would be impossible to connect them into a connected c-graph so that the three e-graphs were really three different e-graphs, and that they were in a genus one u-c-tree.  We analyze how could we possibly connect the ``middle'' e-graph (supposing three parallel  e-graphs) to the other parts of the u-c-tree in order to obtain the desired contradiction. 

Suppose the statement of Theorem 1 was false. Let $\mathcal{E}_1$, $\mathcal{E}_2$, $\mathcal{E}_3$  be three different e-graphs of a genus one u-c-tree $C$ parallel to each other. Let  $\mathcal{P}_1$,  $\mathcal{P}_2$,  $\mathcal{P}_3$ be the  sets of parallel edges, and  $\mathcal{U}_1$, $\mathcal{U}_2$, $\mathcal{U}_3$ be the sets  uncrossed edges of $\mathcal{E}_1$, $\mathcal{E}_2$, $\mathcal{E}_3$, respectively. Let edges $e_1, e_2,\ldots, e_n$ be all of the
edges crossing their parallel edges. 
Let $\mathcal{P}=\mathcal{P}_1\cup \mathcal{P}_2\cup \mathcal{P}_3=\{a_1, \ldots, a_k\}$, where $a_1, \ldots, a_k$ are increasingly parallel. If $a_1=BC$ and $a_k=DA$, then arcs $\widehat{AB}$ and $\widehat{CD}$ are the minimal arcs such that  all edges from $\mathcal{P}$ have
 one of their endpoints on $\widehat{AB}$ while the other on $\widehat{CD}$. Let edge $BC$  be an edge of $\mathcal{P}_1$ and $DA$ an edge of $\mathcal{P}_3$.
Recall that different e-graphs have no common points
and call $\mathcal{E}_1$
the left e-graph,  $\mathcal{E}_3$ the
right e-graph and $\mathcal{E}_2$  the middle e-graph.
Since  e-graphs $\mathcal{E}_1$, $\mathcal{E}_2$ and $\mathcal{E}_3$
 are subtrees of u-c-tree C, they are all connected to each other
within the u-c-tree. We concentrate our efforts on how could  $\mathcal{E}_2$ be connected to the other parts of the u-c-tree $C$.
 \begin{figure}
\begin{center} 
%\epsfbox{le9.eps} 
\epsfbox{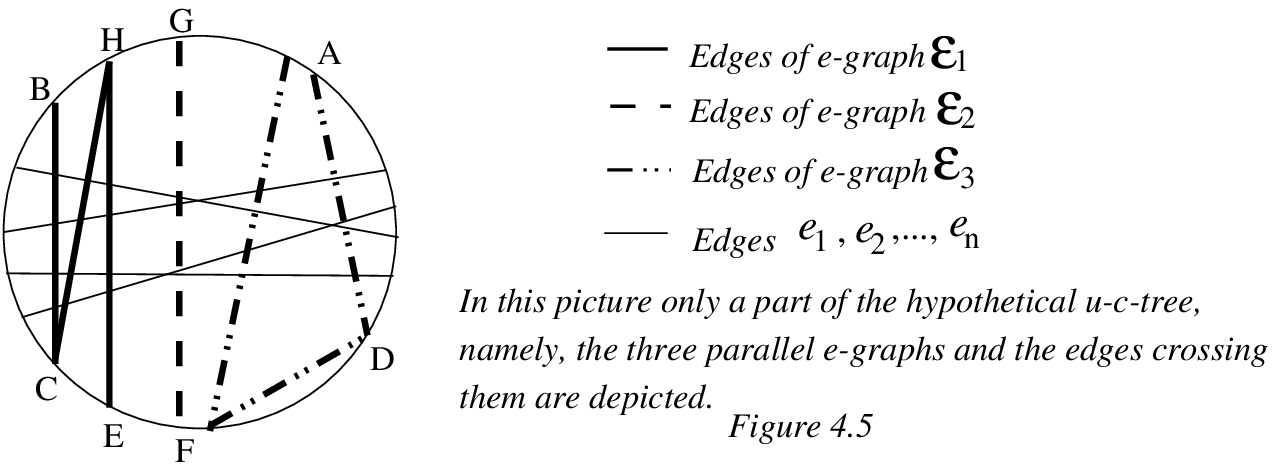}
\end{center} 
\end{figure}

We prove that there is no path connecting   $\mathcal{E}_1$ and  $\mathcal{E}_2$ such that the path contains exclusively edges between   $\mathcal{E}_1$ and  $\mathcal{E}_2$. Analogously, there is no path connecting   $\mathcal{E}_2$ and  $\mathcal{E}_3$ such that the path contains exclusively edges between   $\mathcal{E}_2$ and  $\mathcal{E}_3$.

\begin{figure}
\begin{center}
\epsfbox{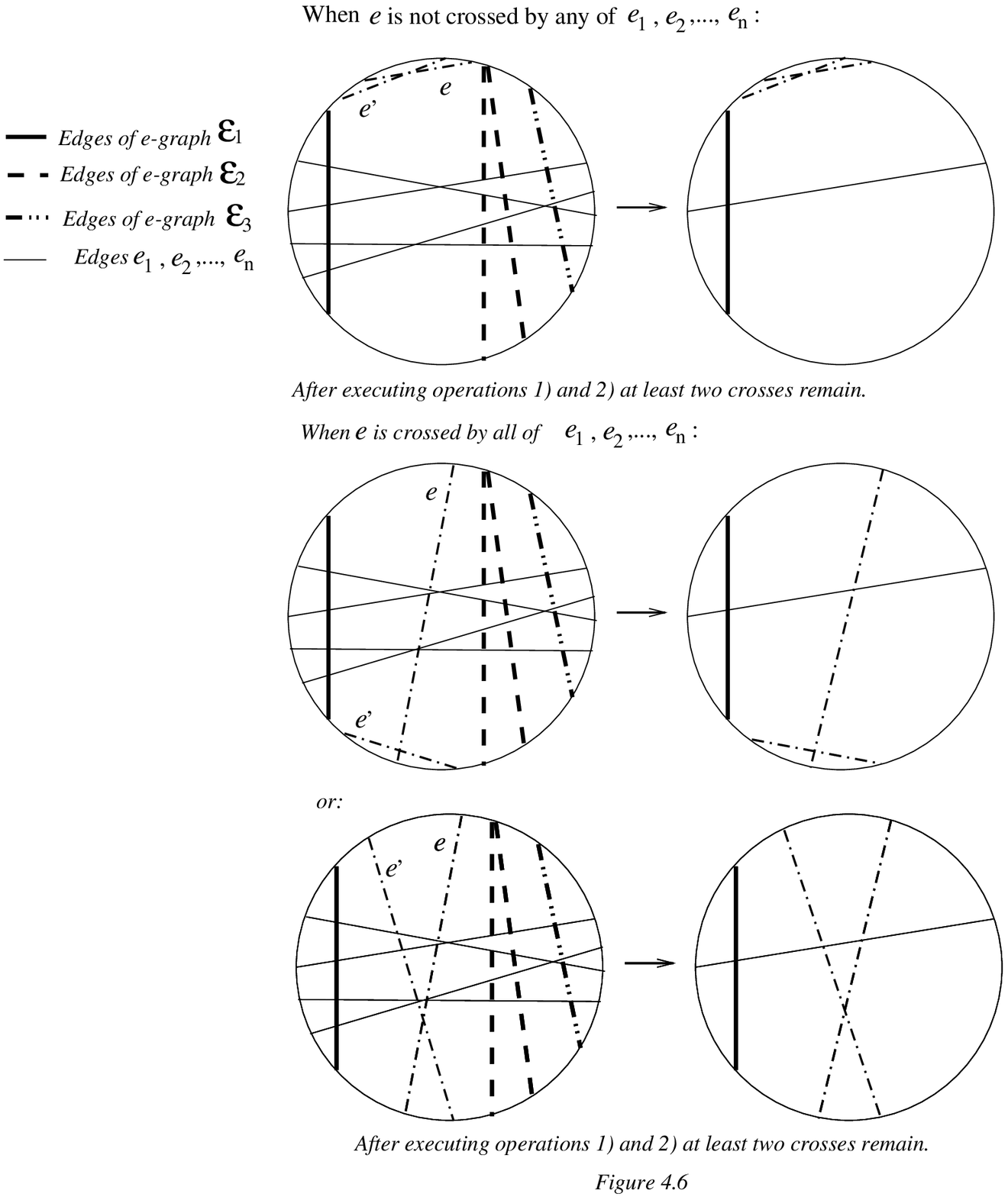} 
 \end{center}
\end{figure}

Suppose the opposite. 
Suppose that there was a path consisting of edges between  $\mathcal{E}_1$ and  $\mathcal{E}_2$ 
 connecting 
$\mathcal{E}_2$  to  $\mathcal{E}_1$. Then, 
 either  $\mathcal{E}_2$  is connected  to $\mathcal{E}_1$
 by only uncrossed edges and edges parallel to  $\mathcal{E}_1$
 or  $\mathcal{E}_2$  is connected  to $\mathcal{E}_1$ 
 with edges among which there are edges which 
 are neither uncrossed nor parallel to $\mathcal{E}_1$. 
 From the definition of  an e-graph we see that only the second possibility 
 might hold.  
 However, if there was an edge $e$ on the path connecting   $\mathcal{E}_1$ and 
   $\mathcal{E}_2$
which is neither uncrossed nor parallel to 
    $\mathcal{E}_1$, then it would have been crossed by some edge which 
     crosses none of $\mathcal{E}_1$, $\mathcal{E}_2$, $\mathcal{E}_3$ or 
it would have been  not crossed by some which crosses these. 
However, if $e$ is not crossed by some of $e_1, e_2, ..., e_n$, 
 then it is not crossed by any of them\footnote{Since if an edge  $e$ between $\mathcal{E}_1$ and $\mathcal{E}_2$ 
 is crossed by some of $e_i$, then $e$'s endpoints are on the two 
 different arcs between the e-graphs, and in the case $e$'s 
 endpoints are on the two different arcs between the e-graphs ($\widehat{EF}$ and $\widehat{GH}$, Figure 4.5), then $e$ 
 is crossed by all $e_1, e_2, ..., e_n$, since the endpoints of 
  $e_1, e_2, ..., e_n$ are left to $\mathcal{E}_1$ and right to 
  $\mathcal{E}_3$ (left and right, referring to Figure 4.5).}, thus, since $e$ is crossed it must be 
 crossed by some edge which cross none of $\mathcal{E}_1$, 
 $\mathcal{E}_2$, $\mathcal{E}_3$. 
Therefore,  edge $e$ is  crossed by some edge $e'$, such that $e' \not\in  \{e_1, e_2, ..., e_n\}$ and 
  $e'$ is  between $\mathcal{E}_1$ and $\mathcal{E}_2$ (so as for the three e-graphs to be parallel, and $e'$ not to be among  $e_1, e_2, ..., e_n$). All possibilities are depicted on Figure 4.6, where we did not present the whole c-tree, but only the edges that are of importance for our proof. 
 Using operations 1) and 2) it would be impossible to obtain Form 1 or Form 2 (the  cross of the three e-graphs and 
 edges $e_1, e_2,..., e_n$ and the other from the crossing of  $e$ and $e'$ would remain).  Therefore, $\mathcal{E}_2$ is not connected to $\mathcal{E}_1$ or $\mathcal{E}_3$ with edges between them.

  Thus, if $\mathcal{E}_2$ was connected to the other parts of the u-c-tree, namely to $\mathcal{E}_1$ and $\mathcal{E}_3$, then there had to be some edge $e$ on the path which connected $\mathcal{E}_2$  to $\mathcal{E}_1$ and 
   $\mathcal{E}_3$ which was not between $\mathcal{E}_1$ and $\mathcal{E}_2$  or $\mathcal{E}_2$  and $\mathcal{E}_3$.
   However, some of such edges $e$ would intersect some of the three e-graphs and not intersect some other, thus this would contradict that the parallel edges of $\mathcal{E}_1$, $\mathcal{E}_2$, $\mathcal{E}_3$ are parallel to each other. 
  Therefore, there cannot be three e-graphs parallel to each other in a genus one u-c-tree.    The statement of Theorem 1 is proven. 
\end{proof}  
 
%\epsfbox{fig12.eps} 

 \section{The E-Reduction Process} 
 
In this section we describe the $\textbf{\textit{e-reduction process}}$ which reduces a u-c-tree C to a \textbf{\textit{reduced}} \textbf{\textit{form}}  which carries enough information of the original u-c-tree C so that from the reduced form of a u-c-tree C we know how many e-graphs C had and how they were connected among each other. 
 
Given a  u-c-tree $T$ perform the following \textbf{\textit{e-reduction}} \textbf{\textit{process}}: 

\textbf{\textit{First Step.}} For all e-graphs of $T$ do the following:
given e-graph $\mathcal{E}$ in $T$ with set $\mathcal{P}$ of parallel edges,
 delete all but one of the edges of $\mathcal{P}$ (operation 2))
 obtaining u-c-graph $T_1$ from $T$.
 
\textbf{\textit{Second Step.}} Delete all the uncrossed edges of $T_1$
 (operation 1))
obtaining u-c-graph $T_2$. (Note that an edge $e$ is uncrossed in $T_1$ if and only if it
was uncrossed in $T$.)

\textbf{\textit{Third Step.}} If in the original u-c-tree T there
was a path consisting of uncrossed edges $E_1E_2$,..., $E_kE_{k+1}$ 
       connecting e-graphs $\mathcal{E}_1$ (with  arcs  $\mathcal{A}_1$, $\mathcal{A}_2$) and $\mathcal{E}_2$ (with arcs
 $\mathcal{B}_1$, $\mathcal{B}_2$)
 in $T$,
      then, if $A_1A_2$ is the edge left from the set of parallel edges of
       $\mathcal{E}_1$ in the \textit{First Step}
such that $A_1 \in \mathcal{A}_1$ and $A_2 \in \mathcal{A}_2$
and if
$B_1B_2$ is the edge left from the set of parallel edges of $\mathcal{E}_2$ in the \textit{First Step}
such that $B_1 \in \mathcal{B}_1$ and $B_2 \in \mathcal{B}_2$
then add edge $A_jB_i$ ($1\le j,i \le 2$) to u-c-graph $T_2$ provided
$E_1$ belongs to arc $\mathcal{A}_j$ and $E_{k+1}$ belongs to $\mathcal{B}_i$.
Do this for all possible paths of uncrossed edges connecting two  e-graphs in $T$. Call the u-c-graph obtained at the end
$T_3$. Note that $T_3$ is a u-c-tree.

To emphasize the importance of u-c-graphs obtained after the \textit{Second} and \textit{Third Step} of the e-reduction process  we give them special names: \textbf{\textit{pre-reduced forms}} and \textbf{\textit{reduced forms}}, respectively, Figure 5.1. The pre-reduced form of a genus one u-c-tree $T$ is the u-c-tree $T_2$ obtained by the execution
of the first two steps of the e-reduction process on $T$.
The reduced form of a genus one u-c-tree $T$ is the u-c-tree $T_3$ obtained by execution
of the e-reduction process on $T$.
We say that u-c-tree $T$ reduces to  u-c-graph $T_3$
if the u-c-tree $T_3$ is  the reduced form of $T$.

%\begin{figure}
\begin{center} 
% \epsfbox{simplificationMIT.eps} 
\epsfbox{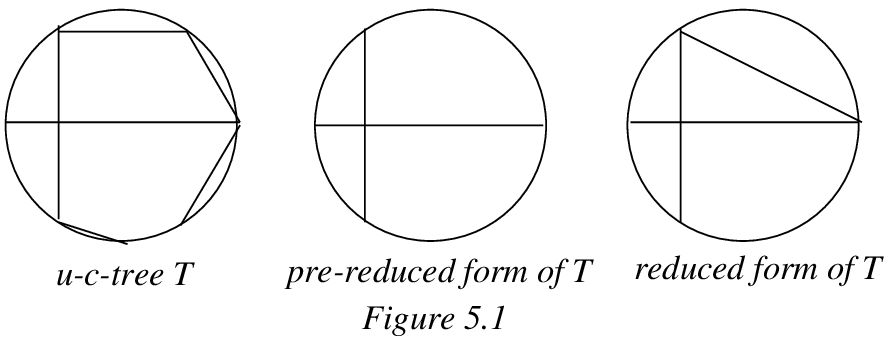}
\end{center} 
%\end{figure}
%\begin{definition}

\begin{definition}
Let $\mathcal{G}$ denote the set of all c-graphs. Let 
$\mathcal{T}_2=\{T_2 \in \mathcal{G} \mid T_2$ is a pre-reduced form of some genus one u-c-tree $C$$\}$ and
$\mathcal{T}_3=\{T_3 \in \mathcal{G} \mid T_3$ is a reduced form of some genus one u-c-tree $C$$\}$

\end{definition}

\begin{lemma}
Any $T_3\in \mathcal{T}_3$ can be obtained from some
 $T_2 \in \mathcal{T}_2$ by addition of uncrossed edges but without addition of any vertices so that the obtained u-c-graph is a tree.
Also, all possible u-c-trees T obtained by taking some u-c-graph 
 $T_2 \in \mathcal{T}_2$ and adding  
exclusively  uncrossed edges without adding any vertices so as to form a tree  are in  $\mathcal{T}_3$.
\end{lemma}
\begin{proof}
Follows from the e-reduction process.
\end{proof}

\begin{lemma}
 Let forms $T_2$ and $T_3$ be the pre-reduced and reduced forms of 
  u-c-tree $T$, respectively. Let  $T_2$ have $n_p$ points and
$n_e$ edges. Then $T_2$ is a subgraph of the tree $T_3$, the set of vertices of $T_2$ is equal to the set of vertices of $T_3$, and there are exactly 
$n_p-n_e-1$ edges of $T_3$ which are not in $T_2$. All of these edges are uncrossed.
\end{lemma}

\begin{proof}
The facts that $T_2$ is a subgraph of $T_3$, $T_3$ is a tree,  the set of vertices of $T_2$ is equal to the set of vertices of $T_3$, and that the edges which belong to $T_3$ but do not belong to $T_2$ are all uncrossed follow directly from the e-reduction  process. Since $T_3$ is a tree on $n_p$ points it has $n_p -1 $ edges, and since $T_2$ has $n_e$ edges there are exactly    $n_p-n_e-1$ in  $T_3$ which are not in $T_2$. 
\end{proof}

\begin{lemma}
Let $T_2$ be the pre-reduced form of a genus one u-c-tree $C$. Then $T_2$ is an offspring of $C$.
\end {lemma}

\begin{proof}
Since the \textit{First Step} and the \textit{Second Step} of the e-reduction process require only repeated execution of operations 1) and 2) on a u-c-tree $C$, the u-c-graph obtained after the \textit{Second Step}, $T_2$, is an offspring of $C$ by definition.
\end{proof}

\begin{proposition}
%The genus of the reduced form $T_3$ of a C-tree $T$
% is equal to the genus of $T$.
The genus of the reduced form $T_3$ of a u-c-tree $T$ is one if and only if
 the genus of $T$ is one.
\end{proposition}

\begin{proof}
Let $T_2$ be the pre-reduced form of u-c-tree $T$. Then $T_2$
is an offspring of $T$ by Lemma 9. Note that
 $T_2$ is an offspring of $T_3$, since deleting the uncrossed edges of $T_3$ results in $T_2$. 
Thus,   $T$ and $T_3$ have a common
offspring and so their final offspring is the same. By  Proposition 1' the genus of
 $T$ is one if and only if the genus of $T_3$ is one.
\end{proof}

\begin{lemma}
If $T_2$ is a pre-reduced form of a genus one u-c-tree T, 
then $T_2$ descends to Form 1 or Form 2.
\end{lemma}
\begin{proof}
From Lemma 9 we know that $T_2$ is an offspring of T. Since T is
genus one if and only if its final offspring is Form 1 or Form 2, and the final offspring
depends only on the starting u-c-graph, it is clear that from any offspring,
so in particular from $T_2$,
Form 1 or Form 2 can be obtained by executing operations 1) and 2)
on the offspring provided T was genus one.
\end{proof}

\begin{figure}
\begin{center}
\epsfbox{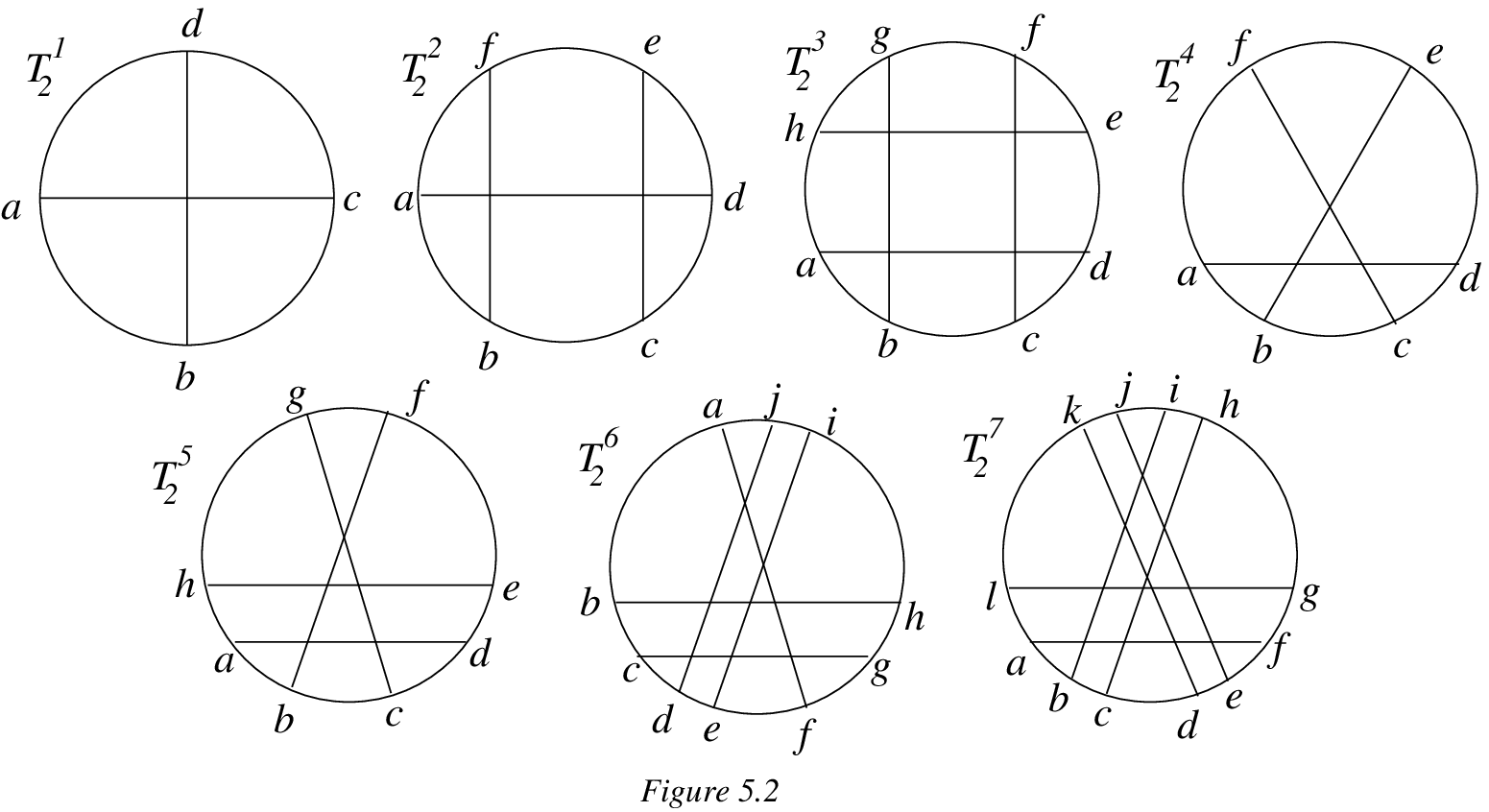}
\end{center}
\end{figure}

\begin{proposition}
If $T_2$ is the pre-reduced form of a genus one u-c-tree T, then $T_2$ 
is one of the 7
u-c-graphs on Figure 5.2\footnote{The labels \textit{a, b, c, d, e, f, g, h, i, j, k, l} are not part of the forms. They serve only to enable us  to specify certain edges of the forms.}

\end{proposition}

\begin{proof}
As a pre-reduced form $T_2$ contains no uncrossed edges, and in u-c-graph $T_2$ each e-graph of $T$ is represented by a single edge,
a representative of the
set of parallel edges of the e-graph in $T$. Since, according to Theorem 1, there can
be at most two parallel e-graphs,   in $T_2$ there
are no three edges parallel to each other. Since no different e-graphs
in $T$ had common vertices (Proposition 3), no edges of $T_2$ have common endpoints.
By Lemma 10 $T_2$ descends to Form 1 or Form 2.
These conditions are fulfilled exclusively  in the presented 7 forms, therefore, $\mathcal{T}_2$ $\subseteq$ $\{$${T_2}^1$, ${T_2}^2$, ${T_2}^3$, ${T_2}^4$,
 ${T_2}^5$, ${T_2}^6$, ${T_2}^7$$\}$ as stated in the proposition.
\end{proof}

\section{
%The Nineteen Possible Reduced Forms of Genus One C-Trees; 
Classification of All Genus One C-Trees} 

In this section we  prove that there exist genus one u-c-trees such that their pre-reduced forms are 
${T_2}^1$, ${T_2}^2$, ${T_2}^3$, ${T_2}^4$,
 ${T_2}^5$, ${T_2}^6$,
but that there is no genus one u-c-tree having  ${T_2}^7$ as its pre-reduced form, that is $\mathcal{T}_2$$=$$\{$${T_2}^1$, ${T_2}^2$, ${T_2}^3$, ${T_2}^4$,
 ${T_2}^5$, ${T_2}^6$$\}$. Throughout the analysis of the seven candidates for pre-reduced forms presented in Proposition 5 we  also describe all of the  possible  reduced
forms of genus one u-c-trees. The interrelations of a u-c-tree and its reduced form  
 enables us to get an insight into the behavior of the number of genus one l-c-trees on $n$ points.

\begin{theorem}
$\mathcal{T}_2$$=\{$${T_2}^1$, ${T_2}^2$, ${T_2}^3$, ${T_2}^4$,
 ${T_2}^5$, ${T_2}^6$$\}$.

$\mathcal{T}_3$$=\{$${T_3}^1[1]$,
${T_3}^2[1]$, ${T_3}^2[2]$, ${T_3}^2[3]$,
${T_3}^3[1]$,
${T_3}^4[1]$, ${T_3}^4[2]$,
${T_3}^5[1]$, ${T_3}^5[2]$, ${T_3}^5[3]$, ${T_3}^5[4]$, ${T_3}^5[5]$, ${T_3}^5[6]$,
${T_3}^6[1]$, ${T_3}^6[2]$, ${T_3}^6[3]$, ${T_3}^6[4]$, ${T_3}^6[5]$, ${T_3}^6[6]$$\}$, Figure 6.1. 
\end{theorem}

\begin{figure}
\begin{center}   
\epsfbox{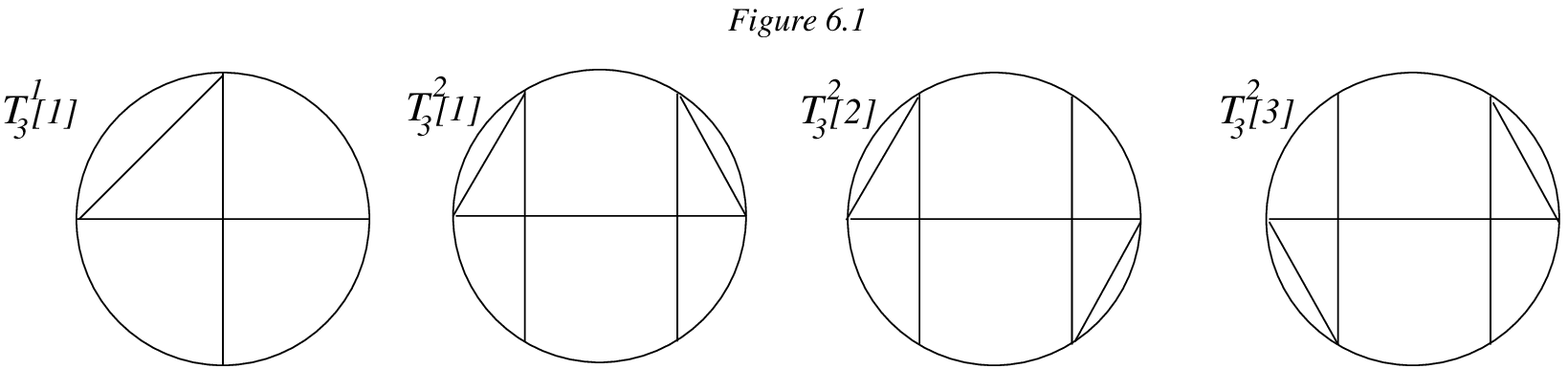} 
\end{center}
\begin{center} 
\epsfbox{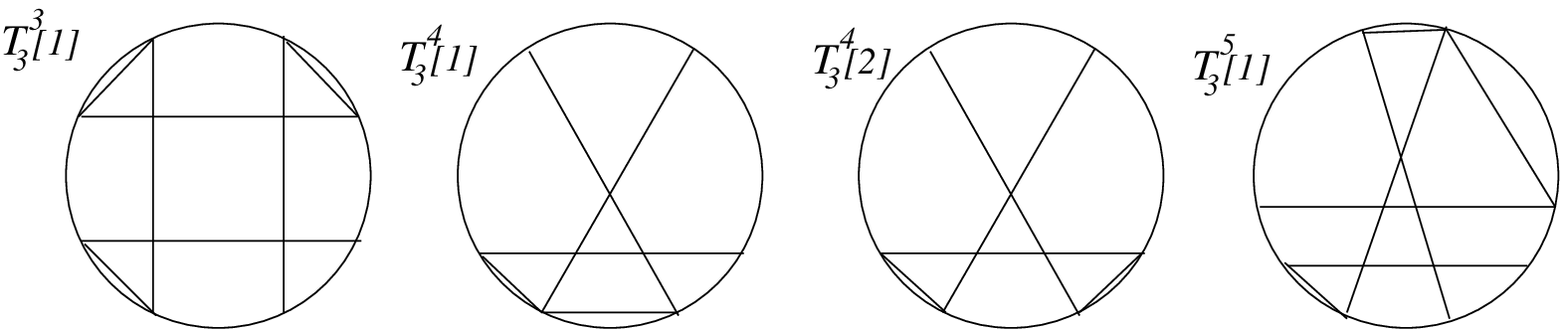} 
\end{center}
 \begin{center} 
\epsfbox{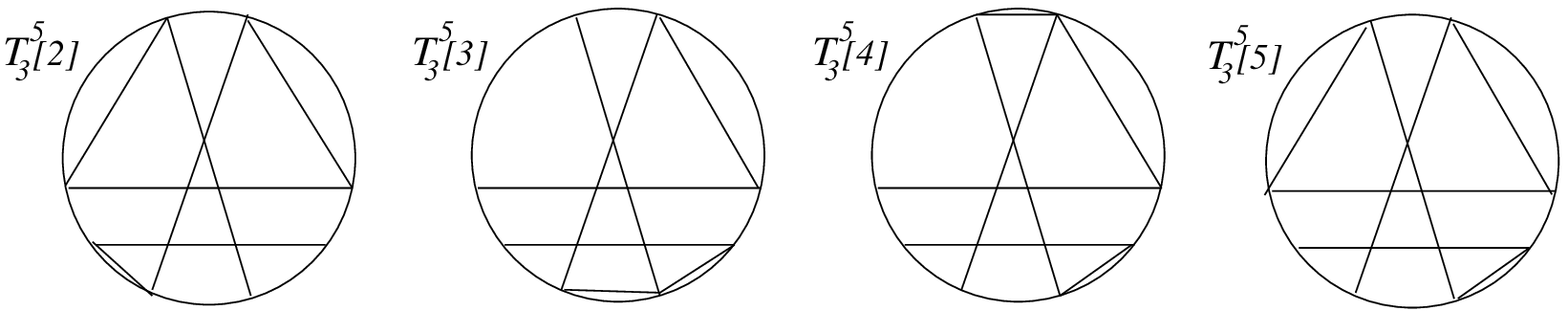} 

\end{center}
\begin{center} 
\epsfbox{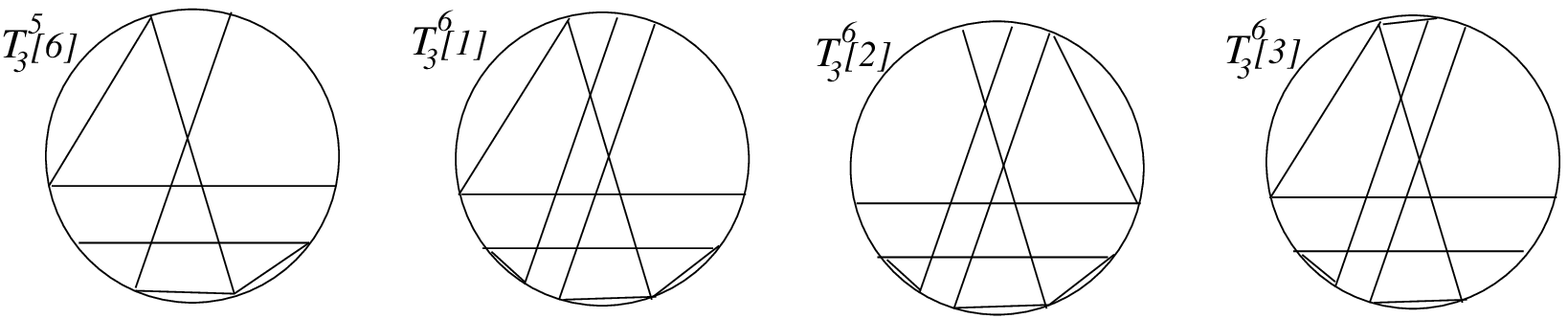} 
 \end{center}
\begin{center} 
\epsfbox{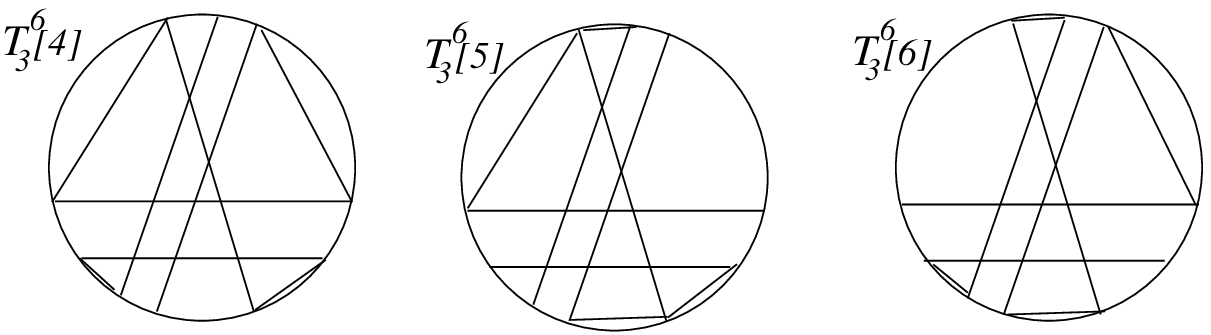}
\end{center} 
%\epsfbox{type6.eps} 
 \end{figure}

%\epsfbox{type7.eps} 
%\epsfbox{type8.eps} 

%\epsfbox{type9.eps} 
%\epsfbox{type10.eps} 

%\epsfbox{type11.eps} 
%\epsfbox{type12.eps} 

%\epsfbox{type13.eps} 
%\epsfbox{type14.eps} 

%\epsfbox{type15.eps} 
%\epsfbox{type16.eps} 

%\epsfbox{type17.eps} 
%\epsfbox{type18.eps} 

%\epsfbox{type19.eps} 

\begin{proof}
We try to obtain a tree from each of 
 the seven different candidates for 
pre-reduced forms from Proposition 5 by  addition of uncrossed edges as described in Lemma 8, pretending the candidates were in fact pre-reduced forms. In case
 we are able to obtain a tree without adding vertices  as Lemma 8 states, we know that the candidate was in fact a pre-reduced form, since the tree we obtain does have the candidate for its pre-reduced form and itself for its reduced form. However, in case we are unable to obtain a tree by the addition of uncrossed edges as  Lemma 8 states, we know that the candidate cannot possibly be a pre-reduced form. Also, for the candidates which prove to be pre-reduced forms we  give all of the reduced forms which can result from them after the execution of the 
\textit{Third Step} of the e-reduction process. Since all reduced forms are obtainable from some pre-reduced form by the execution of the \textit{Third Step} of the e-reduction process we  get all of the reduced forms of genus one u-c-trees.The analysis of the seven candidates follows. We are referring to the figure of Proposition 5.

$\bullet$ ${T_2}^1$ has 4 points and 2 edges. In order to obtain a tree one edge
needs to be added. There are 4 possibilities for uncrossed edges
between the e-graphs: $ab$, $bc$, $cd$ and $da$. 
All the 4 u-c-trees which could be obtained by adding one
of these edges are isomorphic.
Thus, one reduced form: ${T_3}^1[1]$ can be obtained from ${T_2}^1$.

$\bullet$ ${T_2}^2$ has 6 points and 3 edges, thus 2 edges need to be added  to
obtain a tree. There are 4 possibilities for uncrossed edges between
the e-graphs: $ab$, $cd$, $de$ and $fa$. (Note that $bc$ and $ef$ are not such edges, since
$bf$ and $ce$ have to become different e-graphs, and if there would be an uncrossed
 edge connecting these parallel edges, they would form a single e-graph.)
 Edges $fa$ and $ab$ cannot be added at the same time to form  a tree 
 nor can $cd$ and $de$ together,
 \{$fa$, $ab$, $bf$\} or
 \{$cd$, $de$, $ec$\}  would form  a cycle.
 Thus, the 2 edges which could be added are: \{$fa$, $cd$\},
 \{$fa$, $de$\}, \{$ab$, $cd$\}, \{$ab$, $de$\}.
  The  u-c-trees obtained by adding \{$fa$, $de$\} and \{$ab$, $cd$\}
   are isomorphic.
  However, none of the u-c-trees obtained by adding  \{$fa$, $de$\}, \{$fa$, $cd$\},
  \{$ab$, $de$\} are isomorphic, thus three reduced forms:
${T_3}^2[1]$, ${T_3}^2[2]$, ${T_3}^2[3]$ can be obtained from ${T_2}^2$.

$\bullet$ ${T_2}^3$ has 8 points and 4 edges, thus 3 edges need to be added  to
obtain a tree. There are 4 possibilities for uncrossed edges between
the e-graphs: $ab$, $cd$, $ef$ and $gh$.
All the four u-c-trees which could be obtained by adding three
of these edges are isomorphic.
Thus, one reduced form: ${T_3}^3[1]$ can be obtained from ${T_2}^3$.

$\bullet$ ${T_2}^4$ has 6 points and 3 edges, thus 2 edges need to be added to
obtain a tree. There are 6 possibilities for uncrossed edges between
the e-graphs: $ab$, $bc$, $cd$, $de$, $ef$, and $fa$.  
 Edges \{$ab$, $ed$\}, \{$bc$, $ef$\}, \{$cd$, $fa$\}
  cannot added
 at the same time to form a tree,
 since 
 \{$ab$, $be$, $ed$, $da$\},  \{$bc$, $cf$, $ef$, $eb$\} or
\{$cd$, $da$, $fa$, $fc$\} would form a cycle.
 Thus, the 2 edges which might be added are:
 \{$ab$, $bc$\}, \{$ab$, $cd$\}, \{$ab$, $ef$\}, \{$ab$, $fa$\},
 \{$bc$, $cd$\}, \{$bc$, $de$\}, \{$bc$, $fa$\},
 \{$cd$, $de$\}, \{$cd$, $ef$\}, 
 \{$de$, $ef$\}, \{$de$, $fa$\}, \{$ef$, $fa$\}.
 % nahat.. 6 alatt 2  -5 az 12...?!
   The  u-c-trees obtained by adding
 \{$ab$, $bc$\},  \{$ab$, $fa$\}, \{$bc$, $cd$\}, \{$cd$, $de$\},  \{$de$, $ef$\},  \{$ef$, $fa$\} are isomorphic; also  
 u-c-trees obtained by adding
 \{$ab$, $cd$\}, \{$ab$, $ef$\},  \{$bc$, $de$\}, \{$bc$, $fa$\}, \{$cd$, $ef$\}, \{$de$, $fa$\} are
 isomorphic. 
  However,  u-c-trees obtained by adding \{$ab$, $bc$\}, \{$ab$, $cd$\} 
are not isomorphic, thus two reduced forms:
${T_3}^4[1]$, ${T_3}^4[2]$ can be obtained from ${T_2}^4$.
 
$\bullet$ ${T_2}^5$ has 8 points and 4 edges, thus 3 edges need to be added to
obtain a tree. There are 6 possibilities for uncrossed edges between
the e-graphs: $ab$, $bc$, $cd$, $ef$, $fg$, and $gh$. 
Since we want to obtain a tree,  one of \{$ef$, $gh$\} and one of
\{$ab$, $cd$\} must be among the added uncrossed edges in order for 
edges $he$ and $ad$ to be connected to something. Thus, two out of the three
edges needed to be added must be \{$ef$, $ab$\} or
\{$ef$, $cd$\} or \{$gh$, $ab$\}
or \{$gh$, $cd$\}.
In order to obtain a tree the three added edges can be:
\{$ef$, $ab$, $fg$\}, \{$ef$, $ab$, $gh$\}, \{$ef$, $ab$, $bc$\}, \{$ef$, $ab$, $cd$\},
\{$ef$, $cd$, $bc$\}, \{$ef$, $cd$, $fg$\}, \{$ef$, $cd$, $gh$\},
\{$gh$, $ab$, $bc$\}, \{$gh$, $ab$, $cd$\}, \{$gh$, $ab$, $fg$\},
\{$gh$, $cd$, $bc$\}, \{$gh$, $cd$, $fg$\}. 
    The    u-c-trees obtained by adding
\{$ef$, $ab$, $fg$\}, \{$ef$, $ab$, $bc$\} are isomorphic; also  u-c-trees obtained by adding
\{$ef$, $ab$, $gh$\}, \{$ef$, $ab$, $cd$\} are isomorphic; also  u-c-trees obtained by adding
\{$ef$, $cd$, $bc$\}, \{$gh$, $ab$, $fg$\} are isomorphic; also  u-c-trees obtained by adding
\{$ef$, $cd$, $fg$\}, \{$gh$, $ab$, $bc$\} are isomorphic; also  u-c-trees obtained by adding
\{$ef$, $cd$, $gh$\}, \{$gh$, $ab$, $cd$\} are isomorphic; also  u-c-trees obtained by adding
\{$gh$, $cd$, $bc$\}, \{$gh$, $cd$, $fg$\} are isomorphic. 
  However,  u-c-trees obtained by adding \{$ef$, $ab$, $fg$\}, \{$ef$, $ab$, $gh$\},
 \{$ef$, $cd$, $bc$\}, \{$ef$, $cd$, $fg$\}, \{$ef$, $cd$, $gh$\}, \{$gh$, $cd$, $bc$\}
are not isomorphic, thus six reduced forms:
${T_3}^5[1]$, ${T_3}^5[2]$,${T_3}^5[3]$, ${T_3}^5[4]$, ${T_3}^5[5]$,
${T_3}^5[6]$ can be obtained from ${T_2}^5$.

$\bullet$ ${T_2}^6$ has 10 points and 5 edges, thus 4 edges need to be added to
obtain a tree. There are 6 possibilities for uncrossed edges between
the e-graphs: $ab$, $cd$, $ef$, $fg$, $hi$, and $ja$.
There are 15 possibilities to choose 4 edges out of these 6:
             \{$ab$, $cd$, $ef$, $fg$\},
\{$ab$, $cd$, $ef$, $hi$\},
\{$ab$, $cd$, $fg$, $hi$\}, 
\{$ab$, $ef$, $fg$, $hi$\},
\{$cd$, $ef$, $fg$, $hi$\},
\{$ab$, $cd$, $ef$, $ja$\},
\{$ab$, $cd$, $fg$, $ja$\},
\{$ab$, $ef$, $fg$, $ja$\},
\{$cd$, $ef$, $fg$, $ja$\},
\{$ab$, $cd$, $hi$, $ja$\},
\{$ab$, $ef$, $hi$, $ja$\},
\{$cd$, $ef$, $hi$, $ja$\},
\{$ab$, $fg$, $hi$, $ja$\},
\{$cd$, $fg$, $hi$, $ja$\},
\{$ef$, $fg$, $hi$, $ja$\}.
Edges \{$ab$, $bh$, $hi$, $ie$, $ef$, $fa$\} form a cycle thus \{$ab$, $hi$, $ef$\} cannot be in a tree.
Similarly  \{$cd$, $dj$, $ja$, $af$, $fg$, $gc$\} form a cycle thus \{$cd$, $ja$, $fg$\}
cannot be in a tree. Therefore only 9 possibilities out of
15 remain:
\{$ab$, $cd$, $ef$, $fg$\},
\{$ab$, $cd$, $fg$, $hi$\}, 
\{$cd$, $ef$, $fg$, $hi$\},
\{$ab$, $cd$, $ef$, $ja$\},
\{$ab$, $ef$, $fg$, $ja$\},
\{$ab$, $cd$, $hi$, $ja$\},
\{$cd$, $ef$, $hi$, $ja$\},
\{$ab$, $fg$, $hi$, $ja$\},
\{$ef$, $fg$, $hi$, $ja$\}.
 Note that the   u-c-trees obtained by adding
\{$ab$, $cd$, $ef$, $fg$\}, \{$ab$, $fg$, $hi$, $ja$\}
 are isomorphic; also  u-c-trees obtained by adding
\{$cd$, $ef$, $fg$, $hi$\}, \{$ab$, $cd$, $hi$, $ja$\}
are isomorphic; also  u-c-trees obtained by adding 
\{$ab$, $cd$, $ef$, $ja$\}, \{$ef$, $fg$, $hi$, $ja$\} are isomorphic. 
However,  u-c-trees obtained by adding 
\{$ab$, $cd$, $ef$, $fg$\},
\{$cd$, $ef$, $fg$, $hi$\},
\{$ab$, $cd$, $ef$, $ja$\},
\{$ab$, $cd$, $fg$, $hi$\},
\{$ab$, $ef$, $fg$, $ja$\},
\{$cd$, $ef$, $hi$, $ja$\}
are not isomorphic, thus six reduced forms:
${T_3}^6[1]$, ${T_3}^6[2]$,${T_3}^6[3]$, ${T_3}^6[4]$, ${T_3}^6[5]$,
${T_3}^6[6]$ can be obtained from ${T_2}^6$.

$\bullet$ Finally, ${T_2}^7$ has 12 points and 6 edges, thus
  5 edges need to be added to
obtain a tree.
There are 6 possibilities for uncrossed edges between
the e-graphs: $ab$, $cd$, $ef$, $gh$, $ij$, and $kl$. 
Note, however, that
\{$ab$, $bi$, $ij$, $je$, $ef$, $fa$\}, \{$cd$, $dk$, $kl$, $lg$, $gh$, $hc$\} are cycles, thus
\{$ab$, $ij$, $ef$\}, \{$cd$, $kl$, $gh$\} cannot be added to obatin a tree. However,  
there is no way to choose 5 edges out of the 6 possible so that none of the
 sets \{$ab$, $ij$, $ef$\}, \{$cd$, $kl$, $gh$\} is a subset of the set of the chosen 5 edges.
 Therefore, no reduced form can be obtained from  ${T_2}^7$.

We have analyised all possible forms from Proposition 7, and saw that some element  of $\mathcal{T}_3$ can be obtained from  all of ${T_2}^1$, ${T_2}^2$, ${T_2}^3$, ${T_2}^4$,
 ${T_2}^5$, ${T_2}^6$, but no element of $\mathcal{T}_3$ can be obtained from  ${T_2}^7$. Using the result of Proposition 7 we conclude that 
$\mathcal{T}_2$$=\{$${T_2}^1$, ${T_2}^2$, ${T_2}^3$, ${T_2}^4$,
 ${T_2}^5$, ${T_2}^6$$\}$, and since all elements of $\mathcal{T}_3$ are obtainable from some pre-reduced form we have that
$\mathcal{T}_3$$=\{$${T_3}^1[1]$,
${T_3}^2[1]$, ${T_3}^2[2]$, ${T_3}^2[3]$,
${T_3}^3[1]$,
${T_3}^4[1]$, ${T_3}^4[2]$,
${T_3}^5[1]$, ${T_3}^5[2]$, ${T_3}^5[3]$, ${T_3}^5[4]$, ${T_3}^5[5]$, ${T_3}^5[6]$,
${T_3}^6[1]$, ${T_3}^6[2]$, ${T_3}^6[3]$, ${T_3}^6[4]$, ${T_3}^6[5]$, ${T_3}^6[6]$$\}$. This concludes the proof of 
 Theorem 2.

\end{proof}
 
The nineteen  reduced forms from Theorem 2 \textbf{\textit{classify  genus one l-c-trees}}, namely, 
 for all such l-c-trees the corresponding u-c-trees reduce to one 
  of these nineteen  reduced forms.

\section{The Connection Between L-C-Trees and U-C-Trees}

Given a genus one u-c-tree $C$ let \textbf{\textit{l(C)}} denote the number of non-isomorphic genus one l-c-trees corresponding to $C$. 
Alternatively, $l(C)$ is the number of  non-isomorphic l-c-trees obtained from different labelings of 
$C$.  
If  $C_1$, $C_2$, $C_3$, $\ldots$, $C_k$ are all of the  non-isomorphic
genus one u-c-trees on $n$
points, 
  then 
$f(n)=l(C_1)+l(C_2)+l(C_3)+ \cdots+l(C_k)$.

\begin{definition}
Given a l-c-tree $T$ on $n$ points let the rotation of $T$ result in $r(T)$, 
the l-c-tree with edges \{r(a)r(b) $\mid$ $ab$ is an edge of $T$\}, where
$r(i)=i+1$ for $i\in [n-1]$, and $r(n)=1$.
\end{definition}

It directly follows that if $ T$ is a l-c-tree on $n$ points, then $T$ and $r^n(T)$ are isomorphic.

\begin{definition}
Given a u-c-tree $C$  on $n$ points with vertices evenly distributed
on the circle let the rotation of  $C$ result in $R(C)$, 
the u-c-tree obtained by a geometrical rotation of $C$ around the center
of the circle by $\frac{2 \pi}{n}$ in a counterclockwise direction, Figure 7.1.
 \end{definition}
    
%\begin{figure}
\begin{center} 
 \epsfbox{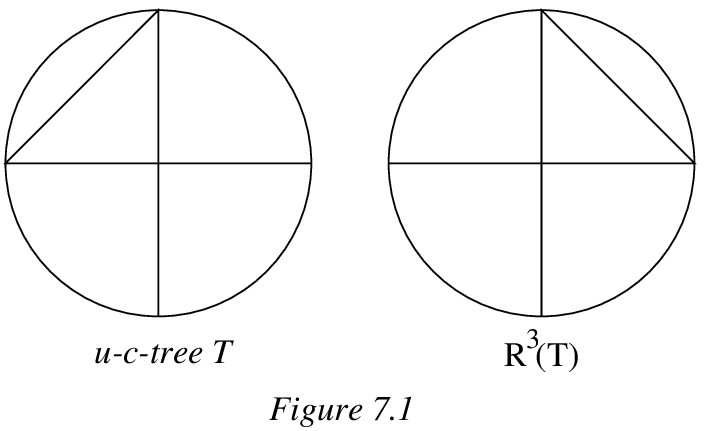} 
\end{center} 
%\end{figure}

Clearly, if $C$ is a u-c-tree on $n$ points with evenly distributed vertices, then $C$ coincides with $R^n(C)$ geometrically.

\begin{proposition}
Given a l-c-tree T on n points, let C be a u-c-tree corresponding to T, such that
C has its vertices evenly distributed on a circle.
%MEGMAGYARAZNI...
Then, $T=r^m(T)$\footnote{When writing $T=r^m(T)$ where $T$ is a l-c-tree
 it is meant that $T$ and $r^m(T)$ are isomorphic.}
if and only if $C=R^m(C)$\footnote{When writing $C=R^m(C)$ where $C$ is a u-c-tree
 it is meant that $C$ and $R^m(C)$ coincide geometrically.}, for all $m\in \mathbb{N}$.
\end{proposition}

%\begin{figure}
\begin{center}
\epsfbox{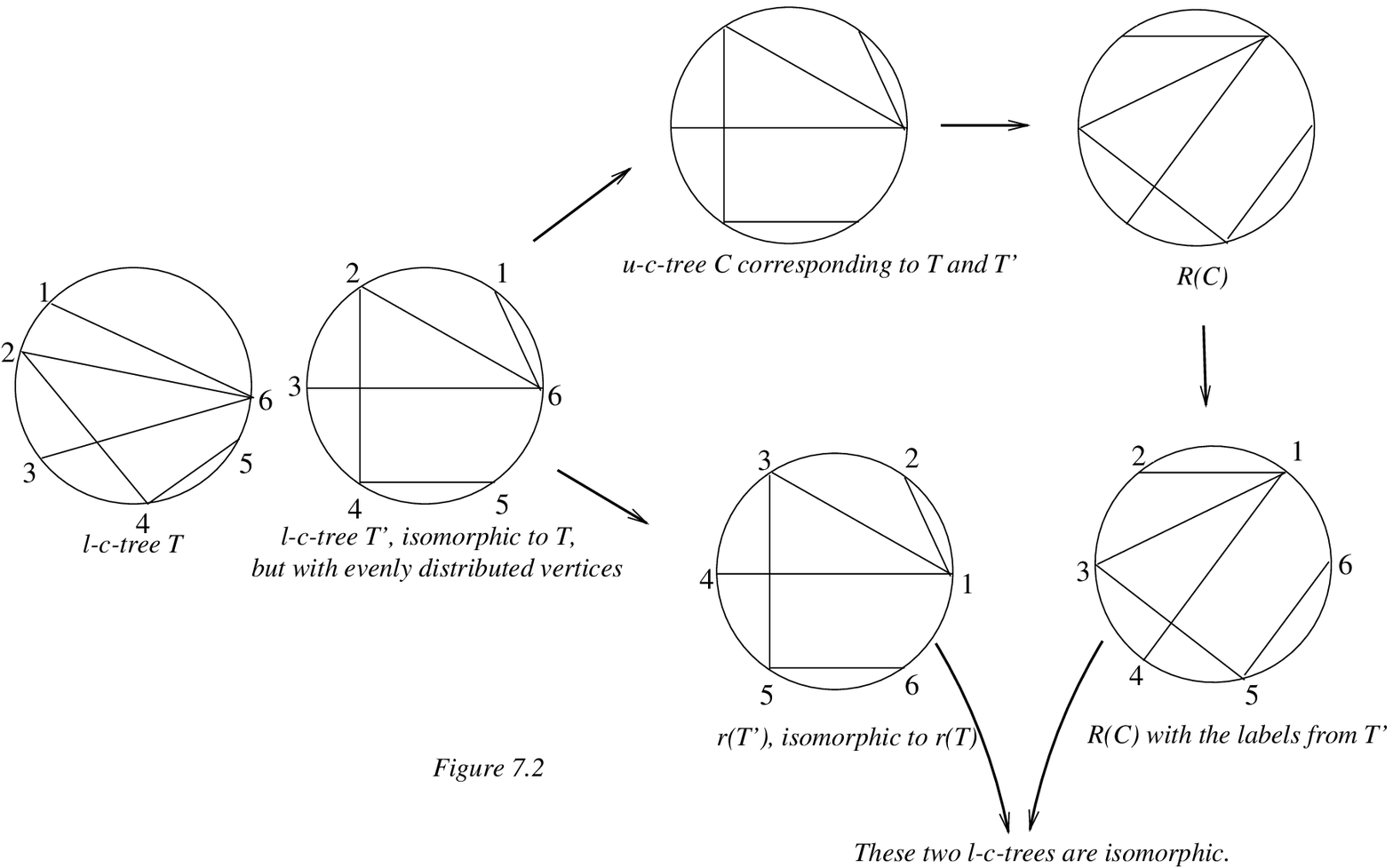}
\end{center}
%\end{figure}

\begin{proof}
Let $T'$ be an  l-c-tree isomorphic to $T$ such that simply deleting the labels of
$T'$ results in $C$. Note that if we fix the labels of $T'$, rotate $C$
and associate the fixed labels to $R(C)$ we obtain the l-c-tree $r(T')$.
(Figure 7.2.)
Iterating this, it is true that associating the fixed labels of $T'$ to
$R^m(C)$ we obtain $r^m(T')$. Thus, $T'=r^m(T')$ if and only if $C=R^m(C)$.
Since $T$ and $T'$ are isomorphic l-c-trees, 
$T=r^m(T)$ if and only if $T'=r^m(T')$. Therefore, $T=r^m(T)$ if and only if $C=R^m(C)$.  
\end{proof}

\begin{proposition}
Given a l-c-tree T on n points let m be the minimal positive integer
such that $T=r^m(T)$. Let C be a u-c-tree corresponding to T.
Then l(C)=m and m divides n.
\end{proposition}

\begin{proof}
Note that $T, r(T), r^2(T),\ldots, r^{m-1}(T)$ are all non-isomorphic l-c-trees,
since if $r^i(T)$ was isomorphic to $r^j(T)$, $0\le i<j\le m-1$ then
$T=r^{j-i}(T)$, $j-i<m$, contradicting that $m$ is the minimal positive
integer such that $T=r^m(T)$.
Also note that any $r^k(T)$, where $k>m-1$, is 
isomorphic to
$r^{k'}(T)$, where $k'$ is the smallest nonnegative remainder of
 $k$ modulo $m$.
Therefore, $T, r(T), r^2(T),\ldots, r^{m-1}(T)$ are all of the
non-isomorphic l-c-trees corresponding to $C$, thus $l(C)=m$.
Since $T=r^n(T)$,  $m$ divides $n$.
(This follows since there exists $k$ such that $0 \le k\cdot m-n < m$, and since $m$ was chosen to be the minimal positive integer with property $T=r^m(T)$, 
and $T=r^{ k\cdot m-n}(T)$ it must be that  $k\cdot m-n=0$.) 

\end{proof}

\begin{proposition}
f(n) is divisible by n if n is prime, n$>$3.
\end{proposition}

\begin{proof}
Let $C_1, C_2,\ldots, C_k$ be all  non-isomorphic u-c-trees on $n$ points.
Let $T_i$ be one of the genus one l-c-trees corresponding to  $C_i$ for
$1\le i\le x$, and let $m_i$ be the minimal positive integer such that
$T_i=r^{m_i}(T_i)$. According to  Proposition 7 all $m_i$ divide $n$, but
it is clear that $m_i$ cannot be 1. Thus $m_i=n$ for all $1\le i\le k$,
and 
$f(n)=l(C_1)+l(C_2)+l(C_3)+ \cdots+l(C_k)=k\cdot n$, is divisible by $n$.
\end{proof}

If for all genus one u-c-trees C on $n$ points $l(C)$  was $n$, it would be 
  trivial to conclude that $f(n)$  is 
  divisible by $n$. However, this is not the case, as already pointed out in Section 2, Figure 2.1.
We  investigate the minimal number of rotations needed for a
u-c-tree $C$ to rotate into itself, which number is equal to $l(C)$ by Proposition 7.
  
%Note that if a u-c-tree $C$ on $n$ points is rotated into itself, then
%its reduced form  has to be rotated into itself as well.

\mbox{}

\noindent \textbf{Observation.}\textit{ When a u-c-tree $C$ is rotated into itself, then its e-graphs
rotate into e-graphs, and paths of uncrossed edges connecting e-graphs
rotate into paths of uncrossed edges connecting e-graphs}.

\mbox{}

Recall that we can think of the reduced form of a u-c-tree C as of a u-c-tree where
e-graphs are represented by e-graphs which have a single edge, and paths of uncrossed edges connecting e-graphs
by single uncrossed edges. 
From the Observation we conclude that if the minimum number of rotations $m$ 
for which a u-c-tree $C$ rotates into itself, 
        $R^m(C)=C$, is less than $n$, and 
        $k$ is the number of vertices of its reduced form, 
        then the reduced form rotates into itself in less then 
        $k$ rotations 
(the points of the reduced form can be made to be 
        evenly distributed on the circle, and the definition of 
        rotation is analogous as in case of a u-c-tree). 
 
\begin{proposition}
 The number of genus one l-c-trees   on $n$ points
 which have corresponding u-c-trees that reduce to reduced forms 
${T_3}^1[1]$,
${T_3}^2[1]$, 
${T_3}^3[1]$,
${T_3}^4[1]$, ${T_3}^4[2]$,
${T_3}^5[1]$, ${T_3}^5[2]$, ${T_3}^5[3]$, ${T_3}^5[4]$, ${T_3}^5[5]$, ${T_3}^5[6]$,
${T_3}^6[1]$, ${T_3}^6[2]$, ${T_3}^6[3]$
 is divisible by $n$.
\end{proposition}

\begin{proof} 
 Let \textbf{\textit{min(T)}} denote the minimum number of rotations needed for a u-c-graph $T$ to rotate into itself.  Since $min(T)<$(the number of vertices of $T$) fails for all reduced forms except for ${T_3}^2[2]$, ${T_3}^2[3]$,  ${T_3}^6[4]$, ${T_3}^6[5]$, and ${T_3}^6[6]$,  only the genus one u-c-trees on $n$ points  having these reduced 
  forms might rotate into themselves in less then $n$ rotations (by  Observation).  
  Therefore, the number of genus one l-c-trees   on $n$ points
 which have corresponding u-c-trees that reduce to
${T_3}^1[1]$,
${T_3}^2[1]$, 
${T_3}^3[1]$,
${T_3}^4[1]$, ${T_3}^4[2]$,
${T_3}^5[1]$, ${T_3}^5[2]$, ${T_3}^5[3]$, ${T_3}^5[4]$, ${T_3}^5[5]$, ${T_3}^5[6]$,
${T_3}^6[1]$, ${T_3}^6[2]$, ${T_3}^6[3]$
 is divisible by $n$, since  for each such
  u-c-tree $C$, $l(C)=n$. 
 \end{proof}
  
When  ${T_3}^2[2]$, ${T_3}^2[3]$,  ${T_3}^6[4]$, ${T_3}^6[5]$, ${T_3}^6[6]$ are rotated into themselves in
less than $k$ rotations (thinking of $k$ as the number of vertices of a particular reduced form) in each of the reduced forms exactly  one edge (which is an e-graph)
 rotates into itself,  while all the other
 edges are ``paired up,'' meaning that  if edge $e$ rotates into edge $e'$ then
 edge $e'$ rotates into edge $e$. 
Consider a u-c-tree $C$ on $n$ points
which  reduces to one of  ${T_3}^2[2]$, ${T_3}^2[3]$,  ${T_3}^6[4]$, ${T_3}^6[5]$, ${T_3}^6[6]$ and for which $min(C)<n$.

Based on the observation made about these reduced forms, we conclude
that in $C$ one of its e-graphs rotates into itself, while the edges not
belonging to that e-graph pair up in the rotation which takes $C$ into
itself in less than $n$ rotations. 
%(This follows since as already noted above, els rotate into els and paths of uncrossed edges into paths of uncrossed edges.)

Let the e-graph of C which rotates into itself have outermost edges 
$e_1$ and $e_2$. In the rotation $e_1$ rotates into $e_2$ and 
$e_2$ rotates into $e_1$, therefore points $Z_i$, $P_i$, $X_j$, $Y_j$  rotate into points $P_i$, 
 $Z_i$,  $Y_j$,  $X_j$ respectively, where $Z_i, P_i, X_j, Y_j$, $i=1, 2,\ldots, r$ and $j=1, 2,\ldots, k$,  are as shown in Figure 7.3.
Thus,  $n=2\cdot(k+r)$ in case a u-c-tree
rotates into itself in less then $n$ rotations, and we
conclude that if for a u-c-tree C on $n$ points $min(C)<n$, then $min(C)= 
\frac{n}{2}$ rotations.
The results we obtained in this discussion are summarized in Theorem 3:

\begin{figure}
\begin{center}
 \epsfbox{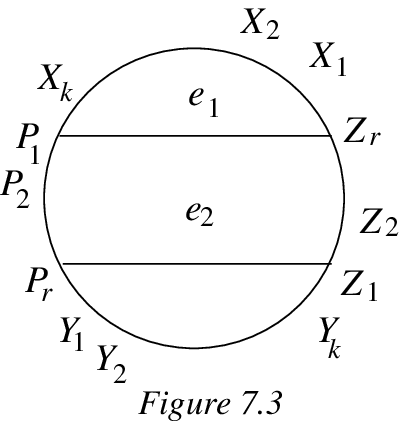} 
\end{center} 
\end{figure}

\begin{theorem}
Given a u-c-tree $C$ on $n$ points, where $n>3$ is odd, $C$ rotates
into itself only after $i\cdot n$ rotations, $i \in \mathbb{N}$, thus $l(C)=n$.
 Given a u-c-tree $C$ on $n$ points, where $n>3$ is even,
if $l(C)\neq n$, then $l(C)=\frac{n}{2}$.
\end{theorem}

\begin{theorem}
 f(n)  is divisible by n for n$>$3, n odd. For n even, n$>$3, f(n) is divisible by $\frac{n}{2}$.
\end{theorem}

\begin{proof}
Theorem 4 follows from Theorem 3, since  $f(n)=l(C_1)+\cdots+l(C_k)$, where $C_1$, $\ldots$, $C_k$ are all of the   non-isomorphic genus one l-c-trees on $n$ points.
\end{proof}

\section{Reduced Forms  ${T_3}^2[2]$, ${T_3}^2[3]$,  ${T_3}^6[4]$, ${T_3}^6[5]$, ${T_3}^6[6]$}

  By Proposition 9 the number of genus one l-c-trees on $n$ points with corresponding u-c-trees that reduce to
${T_3}^1[1]$,
${T_3}^2[1]$, 
${T_3}^3[1]$,
${T_3}^4[1]$, ${T_3}^4[2]$,
${T_3}^5[1]$, ${T_3}^5[2]$, ${T_3}^5[3]$, ${T_3}^5[4]$, ${T_3}^5[5]$, ${T_3}^5[6]$,
${T_3}^6[1]$, ${T_3}^6[2]$, ${T_3}^6[3]$
 is divisible by $n$, thus we are further interested in the 
  number of genus one l-c-trees on $n$ points with corresponding u-c-trees  reducing to
 ${T_3}^2[2]$, ${T_3}^2[3]$,  ${T_3}^6[4]$, ${T_3}^6[5]$, ${T_3}^6[6]$.  
  For each u-c-tree on $n$ points, reducing to  ${T_3}^2[2]$ there is 
  a naturally 
 corresponding u-c-tree on $n$ points reducing to 
${T_3}^2[3]$
 (due to axis symmetry, Figure 8.1) 
 and for each u-c-tree reducing to ${T_3}^6[4]$ there is a 
 u-c-tree reducing to ${T_3}^6[6]$
 for the same reason. This statement is formalized in Proposition 10. 

\begin{center} 
  \epsfbox{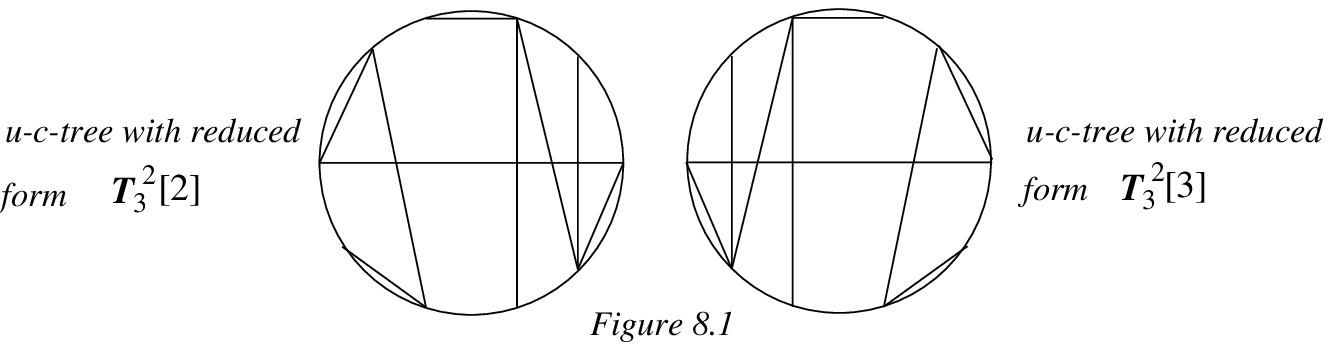} 
\end{center} 

\begin{proposition}

There exists a bijection b  between u-c-trees  on $n$ points reducing to  ${T_3}^2[2]$ (${T_3}^6[4]$)
and u-c-trees  on $n$ points reducing to ${T_3}^2[3]$  (${T_3}^6[6]$) such that if $b(C[1])=C[2]$, then  $min(C[1])=min(C[2])$ to rotate into itself. 

\end{proposition}

\begin{proof}
We exhibit an explicit  bijection $b$ satisfying the conditions of Proposition 10.
Given a u-c-tree $C[1]$  on $n$ points reducing to  ${T_3}^2[2]$ (${T_3}^6[4]$), label its vertices with $\it{1}$ through
$n$ in a counterclockwise direction starting by labeling an arbitrary vertex 
with $\it{1}$.
Take a circle $c_2$ and label $n$ of its points in a clockwise direction with $\it{1}$
through $n$. Draw edges $ij$ on $c_2$ provided some edge of $C[1]$ was
labeled with $ij$. Then delete the labels of the points of $c_2$. The
obtained graph $b(C[1])=C[2]$ is a u-c-tree reducing to ${T_3}^2[3]$  (${T_3}^6[6]$).
Bijection $b$ satisfies the conditions specified. 
\end{proof}

\begin{proposition}
 The number of genus one l-c-trees on $n$ points
  with u-c-trees  reducing to ${T_3}^2[2]$, ${T_3}^2[3]$ ${T_3}^6[4]$, or ${T_3}^6[6]$ is divisible by $n$.
\end{proposition}
\begin{proof}
Proposition 10 proves the existence of a bijection $b$ 
 between u-c-trees  on $n$ points reducing to ${T_3}^2[2]$ (${T_3}^6[4]$)
and u-c-trees  on $n$ points  reducing to ${T_3}^2[3]$  (${T_3}^6[6]$) such that if $b(C[1])=C[2]$, then $min(C[1])=min(C[2])$.
Since for any u-c-tree C, $min(C)=n$ or $min(C)=\frac{n}{2}$, it follows that $l(C[1])+l(C[2])$ is always divisible by $n$ (since $min(C)=l(C)$).
  Therefore, summing $l(C)$ over all u-c-trees $C$ on $n$ points reducing to ${T_3}^2[2]$, ${T_3}^2[3]$ ${T_3}^6[4]$, ${T_3}^6[6]$ we obtain a number divisible by $n$. The statement of the proposition follows. 
\end{proof}

$\bf{Remark.}$ From Proposition 9 and Proposition 11  we conclude that 
  it depends only upon the number of 
 genus one l-c-trees T with corresponding  u-c-trees C  reducing to ${T_3}^6[5]$ and such that $min(C)=\frac{n}{2}$, whether  $f(n)$ is
 divisible by $n$, or only by $\frac{n}{2}$.
 Note that this question is for $n$ even, since for $n$ odd we already saw that $l(C)=n$ for all genus one u-c-trees $C$ on $n$ points.

\section{The Examination of Reduced Form ${T_3}^6[5]$}

 Let $P_n=\{C \in \mathcal{G}$ $\mid$ \textit{C is a u-c-tree on n points reducing to} ${T_3}^6[5]$  \textit{such that} $min(C)=\frac{n}{2}\}$. 
 
\noindent From the Remark of Section 8 $f(n)$ ($n$ even) is divisible by $n$ if and only if 
 $\mid$$P_n$$\mid$ is even, since $f(n)\equiv$$\mid$$P_n$$\mid$$\cdot$$\frac{n}{2}$ (mod $n$). If  $\mid$$P_n$$\mid$ is odd, then $f(n)$  is not divisible by $n$ but is divisible
 by $\frac{n}{2}$.
In this section we investigate the 
 parity of $\mid$$P_n$$\mid$.
\begin{figure}  
\begin{center} 
\epsfbox{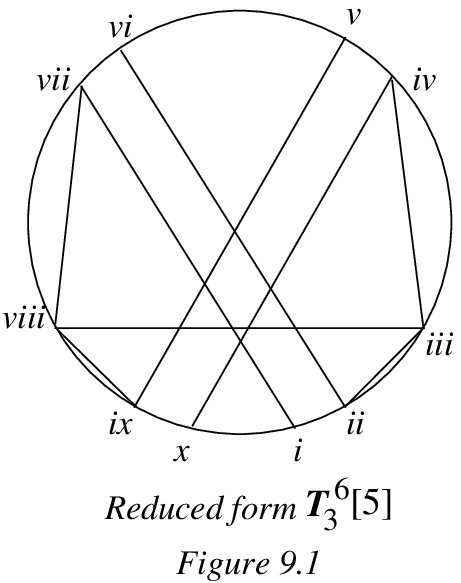} 
\end{center} 
\end{figure}

We say that an \textbf{\textit{e-graph reduces to an edge}} $ab$ if and only if in the e-reduction process edge $ab$ (or $a-b$) was the one left (not deleted) from the set of parallel edges of the e-graph.
In a u-c-tree $C\in$ $P_n$ e-graph $a-b$, $a, b \in \{i, ii, iii, \cdots, ix, x\}$,  is the e-graph
 reducing to edge  $a-b$ in the reduced form of $C$,  Figure 9.1.

 \begin{proposition}
 There exists a bijection $h$ mapping the u-c-trees $C\in$ $P_n$ with
 e-graph $iii-viii$ consisting of more than one edge
 into the u-c-trees $C\in$ $P_n$  
 with e-graph $iii-viii$ consisting of more than one edge, such that if $h(C[1])=C[2]$ then $h(C[2])=C[1]$ and
 $C[1]$ and $C[2]$ are non-isomorphic.
 \end{proposition}

\begin{proof}   
 Given  a u-c-tree $C[1]\in$ $P_n$ we know $min(C)=\frac{n}{2}$ and thus its e-graph $iii-viii$ also 
 rotates into itself in  $\frac{n}{2}$ rotations.
Let $e_1$ and $e_2$ be the e-graph's outermost edges in $C[1]$.
Let  points $X_1, X_2, X_3,..,X_m$ be the points on one arc of the e-graph 
 in  counterclockwise direction
  and let $Y_1,Y_2,Y_3,...,Y_m$ be the points on the other arc of the e-graph 
in  clockwise direction, so $\{e_1, e_2\}=\{X_1Y_1, X_2Y_2\}$.
(Note that the number of points on the two arcs of the e-graph is the same
since they rotate into each other being that $e_1$ rotates into $e_2$
and vica versa. Also, $m>1$ since we supposed the e-graph consists of more than
one edge.).
Leaving the edges of u-c-tree $C[1]$ the same,
except changing the edges of form ($X_i$, $X_j$)
into ($Y_i$, $Y_j$), ($Y_i$, $Y_j$) into
($X_i$, $X_j$) and ($X_i$, $Y_j$) into ($Y_i$, $X_j$) we obtain a u-c-tree
$C[2]$ such that $min(C[2])=\frac{n}{2}$, Figure 9.2. It is clear from the construction that 
if  $h(C[1])=C[2]$ then $h(C[2])=C[1]$. To complete the proof, we need to show that $C[1]$ and $C[2]$ are different. Suppose the opposite. Then e-graph $iii-viii$ also coincides in the two u-c-trees.  Let $e_1=E_{1,1}E_{1,2}$ and $e_2=E_{2,1}E_{2,2}$. Then $e_1$ and $e_2$ need to be connected by the edges of e-graph  $iii-viii$ since they are the outermost edges of the e-graph. Let there be a path of edges of the e-graph connecting vertices $E_{1,i_1}$ and $E_{2, j_1}$ ($1\leq i_1, j_1 \leq 2$) in $C[1]$. Then in $C[2]$ vertices   $E_{1,i_2}$ and $E_{2, j_2}$ ($1\leq i_2, j_2 \leq 2$, $i_1\neq i_2$, $j_1\neq j_2$) would be connected, and if $C[1]=C[2]$ then there would be a circle of edges containing vertices $E_{1,1}$, $E_{1,2}$, $E_{2,1}$, and $E_{2,2}$. This, however, cannot happen since in a tree there cannot be circles.   
\end{proof}

\begin{figure}
  \begin{center}
\epsfbox{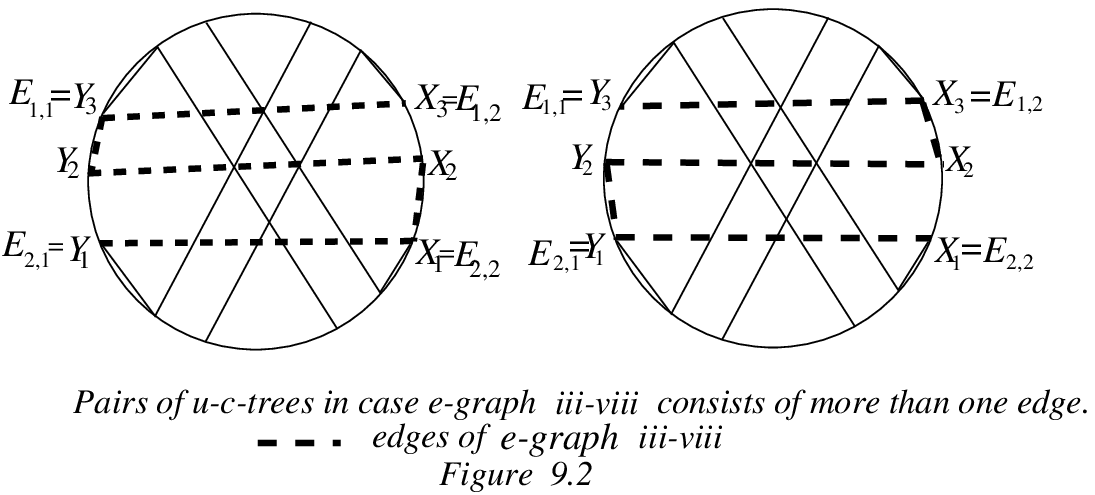} 
\end{center}
\end{figure}

\begin{proposition}
The number of u-c-trees $C$ such that $C\in$ $P_n$ and e-graph $iii-viii$
in $C$ consists of more than one edge is even.
\end{proposition}

\begin{proof}
Follows from the existence of bijection $h$ described in Proposition 12.
\end{proof}

%($de$ ez nem is kell: 
  Proposition 13 implies that
the number of genus one l-c-trees on $n$ points
  with corresponding u-c-trees that are in $P_n$ and for which in the corresponding 
 u-c-trees e-graph $iii-viii$ consists of more than one edge is divisible by $n$.

\mbox{}

\noindent \textbf{ About  U-C-Trees $C\in$ $P_n$ Such That E-Graph $iii-viii$ of $C$
Consists of a Single Edge}

\mbox{}

\noindent It remains to examine the u-c-trees $C\in$ $P_n$ such that e-graph $iii-viii$ of $C$
consists of a single edge. Until the end of this section it is assumed
that all u-c-trees are such.
The definitions assume this as well. Also, we simply denote such a u-c-tree by $C$. The labels $i, ii,\ldots, viii, ix, x$ refer to Figure 9.1 or to a u-c-tree C which has  ${T_3}^6[5]$ as its reduced form.

When we refer to \textbf{\textit{uncrossed edges which connect to e-graph}} $\mathcal{E}$, we mean the set of uncrossed edges $E$, defined recursively as follows:

$\bullet$ all the uncrossed edges which are not  edges of $\mathcal{E}$ but have some of the endpoints of the outermost edges of $\mathcal{E}$ as their endpoints are in $E$

$\bullet$ edge $e$ is in $E$ if it is not in $\mathcal{E}$ and if it is uncrossed and  has a common endpoint with some edge already in $E$

\begin{definition}
Let $K_1$  be the subgraph of $C$ such that it consists of
e-graph $vi-ii$, and  all the uncrossed edges connecting to e-graph $vi-ii$ with the restriction that they have both of their endpoints on arcs $i-iii$ and $v-vii$.
\end{definition}

\begin{definition}
Let $K_2$ be the subgraph of $C$ such that it consists of
e-graph $x-iv$, all of the uncrossed edges connecting to e-graph $x-iv$ with the restriction that they have both of their endpoints on arcs $ix-i$ and $iii-v$.
 \end{definition}
 
\begin{lemma}
E-graph $iii-viii$, $K_1$, and $K_2$  
  uniquely determine a u-c-tree C such that $min(C)=\frac{n}{2}$.
  \end{lemma}
  \begin{proof}
When $C$ rotates into itself in $\frac{n}{2}$ rotations,
   then $K_1$ rotates into the union of e-graph $i-vii$ and the uncrossed edges  
connecting to the e-graph $i-vii$ with the restriction that they have both of their endpoints on arcs $vi-viii$ and $x-ii$,  and $K_2$  
   rotates into the union of e-graph $v-ix$, 
the uncrossed edges connecting to the e-graph $v-ix$ with the restriction that they have both of their endpoints on arcs $viii-x$ and $iv-vi$. Therefore, e-graph $iii-viii$, $K_1$, $K_2$,
and the edges into which $K_1$ and $K_2$ rotate constitute all of the edges
of $C$, thus the statement of the lemma follows.
\end{proof}

\begin{definition}
Given $C$ we say that 
 $K_1$ and $K_2$ are identical if when we label
 the points of $C$ with 1 through n in a clockwise direction so that 
 iii  gets labeled with 1 the edges of $K_1$
 get the same labels as the edges of $K_2$
 when we label
 the points of $C$ with 1 through n in a counterclockwise direction so that
  iii gets
  labeled with 1.
 If this is not the case we say that   $K_1$ and $K_2$ are different.
\end{definition}

Intuitively, $K_1$ and $K_2$ are different if by reflecting
$K_1$ upon the axis of symmetry parallel to $viii-iii$ 
we do not get $K_2$, Figure 9.3. 
 
\begin{figure}
\begin{center}
\epsfbox{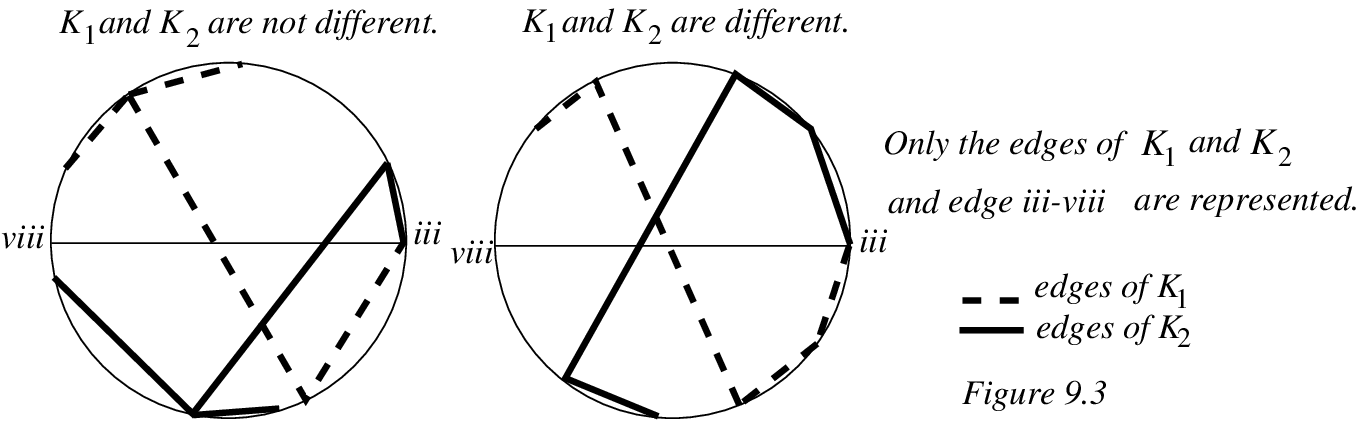} 
 \end{center}
\end{figure}

\begin{lemma}
If in a u-c-tree $C$ $K_1$ is different from $K_2$, then there
exists a bijection $g$
mapping  the set of u-c-trees $C$ with $K_1$ and $K_2$ different into
itself,
 so that if $g(C[1])=C[2]$, then $g(C[2])=C[1]$ and $C[1]$ is not isomorphic to $C[2]$.
 \end{lemma}
\begin{proof}
 Intuitively, $g$ maps u-c-tree $C[1]$ into a u-c-tree $C[2]$ such that
 $C[2]$ is obtained from $C[1]$ by reflecting all the edges of $C[1]$
 upon edge $iii-viii$, Figure 9.4.

\begin{figure}
\begin{center}
\epsfbox{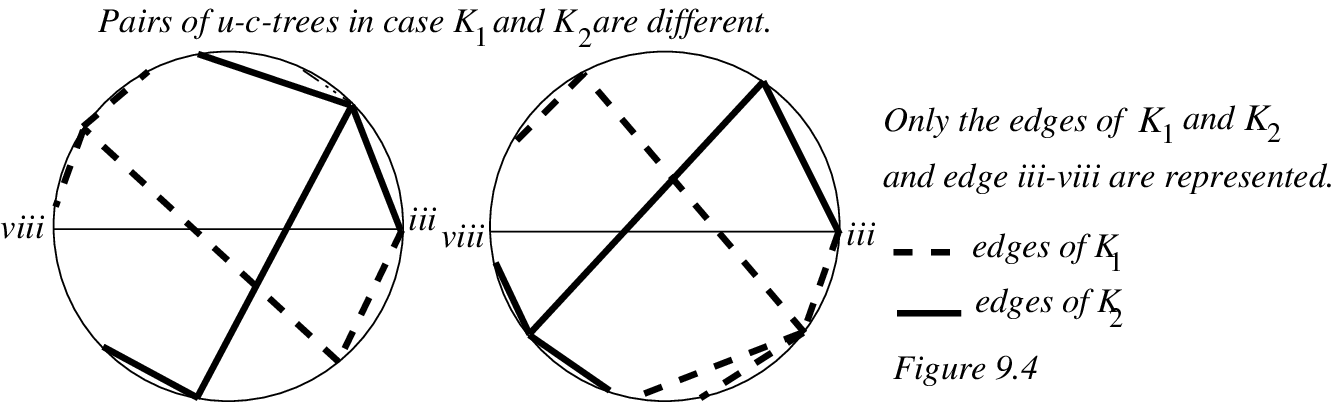}
\end{center}
\end{figure}

 Formally, the following bijection $g$ has the property described.
Given $C[1]$ with $K_1$ and $K_2$ different, label its vertices with $\it{1}$ through
$n$ in a counterclockwise direction starting by labeling $iii$ with $\it{1}$.
Take a circle $c_2$ and label $n$ of its points in a clockwise direction with $\it{1}$
through $n$. Draw edges $ij$ on $c_2$ provided some edge of $C[1]$ was
labeled with $ij$. Then delete the labels of the points of $c_2$. The
obtained graph is $C[2]$.
\end{proof}

\begin{lemma}
The number of u-c-trees C with $K_1$ and $K_2$ different is even.
\end{lemma}
\begin{proof}
The statement follows from the existence of  bijection $g$ described in
Lemma 12.
\end{proof}

     Therefore, it remains to determine the parity  of the number $p$ of 
u-c-trees in $P_n$, such that e-graph $iii-viii$ is a single edge, and 
$K_1$ and $K_2$ are identical.  
 This can happen only if $\frac{n-2}{2}$ is even (since if
$K_1$ and $K_2$ are identical, then the number of edges in $K_1$ and $K_2$
is the same, thus the sum of the number of edges
in $K_1$ and $K_2$ is even, and on the other hand it is  $\frac{n-2}{2}$).
Therefore, in case
$\frac{n-2}{2}$ is odd (which is equivalent to $n$ divisible by \textit{4})
$p=0$ and it follows that:
\begin{theorem}
If $n$ is divisible by 4, the parity of the number of u-c-trees on n
points is even, and so f(n)  is divisible by n.
\end{theorem}

 For $n\equiv$\textit{ 2} (mod \textit{4}) ($\frac{n-2}{2}$ is divisible by \textit{2}) 
 $p$ is equal to the number of possible ways of constructing
 the union of e-graph $vi-ii$ and  all the uncrossed edges connecting to e-graph $vi-ii$ with the restriction that they have both of their endpoints on arcs $i-iii$ and $v-vii$ using  $\frac{n-2}{4}$ edges, since this uniquely determines a u-c-tree $C$ where $iii-viii$ is a single edge and $K_1$ and $K_2$ are identical. We formalize the previous notion 
in the following definition.

\begin{definition}
Suppose we draw e-graph $\mathcal{E}$ with parallel edges represented by vertical lines, then we call the outermost edge on the right the rightmost edges of $\mathcal{E}$.
$\mathcal{K}_m=\{K \in \mathcal{G} \mid$ K is a 
 u-c-graph with $m$ edges consisting of an e-graph $\mathcal{E}$ and uncrossed edges connecting to the e-graph,
with the special property that there is at least one uncrossed edge AX,
 not an edge of the e-graph $\mathcal{E}$, connecting to the e-graph
  whose  rightmost edge is AB and X is in open arc $\widehat{AB}\}$.  
  \end{definition}

Clearly, e-graph $vi-ii$ of $C$ is imitated by e-graph $\mathcal{E}$, that there are 
uncrossed edges on the 
arc between $ii$ and $iii$ is ensured by the edge $AX$ mentioned in the definition, and the uncrossed edges connecting to e-graph  $\mathcal{E}$ imitate the 
 uncrossed edges connecting to e-graph $vi-ii$ of $C$  with the restriction that they have both of their endpoints on arcs $i-iii$ and $v-vii$, and vice versa. From this  $p=\mid \mathcal{K}_{\frac{n-2}{4}}\mid$.

\begin{lemma}
There is an even number of $K \in \mathcal{K}_m$ such that the e-graph of $K$  consists of more than one edge.
\end{lemma}

\begin{proof} 

Let $K'\in \mathcal{K'}_m$ be given, where $\mathcal{K'}_m=\{K \in \mathcal{K}_m \mid$ \textit{the e-graph of K consists of more than one edge\}}. 
Define a function $k:\mathcal{K'}_m \rightarrow \mathcal{K'}_m$ as follows. 
 Think of $K'$ using $k$ edges as $K_1$ with the property that  e-graph $vi-ii$ has more than one edge.
   Let the points of  e-graph $vi-ii$ be $X_1$, $X_2$, $X_3$, $\ldots$, $X_{m_1}$ counterclockwise on one of its arcs  and $Y_1,$ $Y_2,$ $Y_3,$$\ldots$, $Y_{m_2}$ clockwise on the other  of its arcs
 ($m_1>1$ or $m_2>1$). Changing edges of form ($X_i$, $X_j$) 
  into ($Y_i$, $Y_j$), ($Y_i$, $Y_j$) into ($X_i$, $X_j$) and 
  ($X_i$, $Y_j$) into ($Y_i$, $X_j$), and leaving the other edges of $K'$ the same we    
  obtain  $K''=k(K')$  different from $K'$ (The proof that $K'$ and $K''$ are different follows the lines of an analogous proof at the end of Proposition 12.). Clearly, $K'=k(K'')$.  Therefore, 
$\mid \mathcal{K'}_m \mid$ is even, which is the statement of the lemma.
See Figure 9.5. 
 \end{proof}

\begin{figure}
\begin{center}
\epsfbox{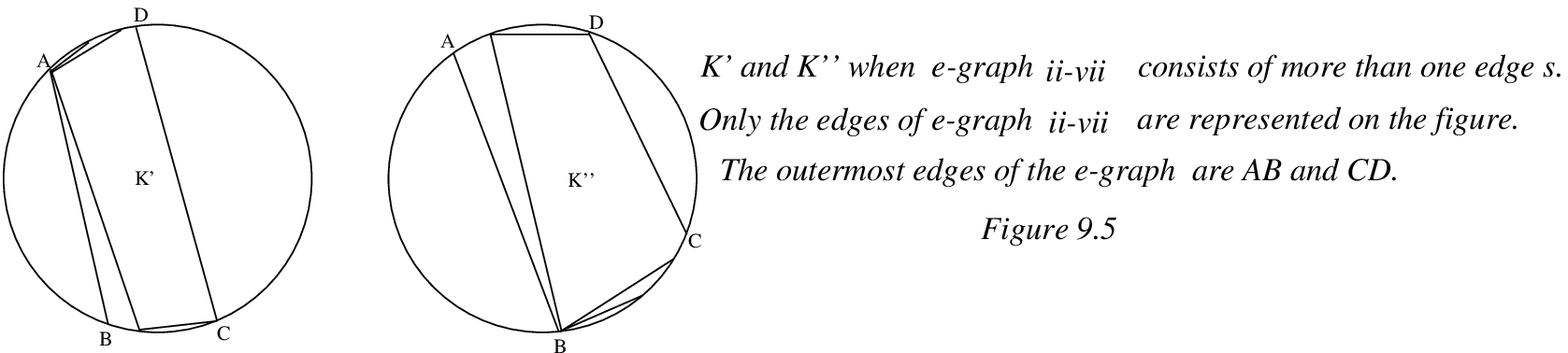} 
 \end{center}
\end{figure}

$\bf{Conclusion.}$ Since $\mid \mathcal{K'}_{\frac{n-2}{4}}\mid$  contributes an even number to $p$,  the only remaining case to consider is when   e-graph $ii-vi$ in $K \in \mathcal{K}_{\frac{n-2}{4}}$ consists of one edge. 
       Our goal is to determine the parity of the number of 
possible constructions of:
 union of e-graph $vi-ii$, consisting of a single edge,  and  all the uncrossed edges connecting to e-graph $vi-ii$ with the restriction that they have both of their endpoints on arcs $i-iii$ and $v-vii$ and that there is at least on edge on arc $ii-iii$, using 
 $\frac{n-2}{4}$ edges.

\section{L-C-Trees and Left-Right Trees}

Observe that the structure described in the Conclusion of the previous section  naturally  decomposes into three parts, namely the edge $vi-ii$, and the uncrossed edges with endpoints on arc $i-iii$ and the uncrossed edges with endpoints on arc $v-vii$. The following two definitions capture the latter two trees, whereas the third definition is a formal way of defining the structure from the Conclusion. We note in advance that the edge $AB$ in the definitions represents the edges $vi-ii$ from the Conclusion.

\begin{definition}
$\mathcal{A}$$=\{T \in \mathcal{G} \mid$ $T$ is a u-c-tree, such that given a specially designated edge $AB$, called the axis of $T$, the union of $T$ and $AB$ is a u-c-tree with all edges uncrossed, and there is no edge in $T$ having $B$ as its endpoint (but clearly, vertex $A$ is a vertex of $T$)$\}$.

If the axis of $T$ is $AB$ and $A$ is a vertex of $T$, then  
vertex  $A$ is called the root of $T$. It is also said that $T$ is rooted at $A$.  (Note, we don't consider $AB$ to be an edge of $T$.)
\end{definition} 

\begin{definition}
$\mathcal{B}$$=\{T \in \mathcal{G} \mid T$ is a u-c-tree, such that given a specially designated edge $AB$, called the axis of $T$, the union of $T$ and $AB$ is a u-c-tree with all edges uncrossed, and there is no edge in $T$ having $A$ as endpoint
 (but clearly, vertex $B$ is a vertex of $T$); furthermore, we require that there is at least one edge of T of form $BC$ with $C$ in open arc $\widehat{BA}\}$. 
If the axis of $T$ is $AB$ and $B$ is a vertex of $T$, then
vertex  $B$ is called the root of $T$. It is also said that $T$ is rooted at $B$. (Note, we don't consider $AB$ to be an edge of $T$.)

%Call any $T \in \mathcal{B}$ an object of type $\mathcal{B}$.
\end{definition} 
 
  \begin{definition}
 $\mathcal{L}=\{L \mid$ $L=(AB$, $T_1$, $T_2)$, where $AB$ is a specially designated edge on the circle, $T_1 \in \mathcal{A}$ rooted at vertex $A$ having axis $AB$, and  $T_2 \in \mathcal{B}$ rooted at vertex $B$ having axis $AB$; the union of  $AB$, $T_1$ and $T_2$ is a u-c-tree with all edges uncrossed, such that $T_1$ and $T_2$ have no vertices in common$\}$.

The number of edges of $L=(AB$, $T_1$, $T_2)\in$ $\mathcal{L}$ is the number of edges of $T_1$ plus the number of edges of $T_2$ plus one.
%Call any $L\in$ $\mathcal{L}$ an object of type  $\mathcal{L}$.

\end{definition}

As already noted at the beginning of the section, the structure defined in the conclusion of Section 9 is modeled by $L=(AB$, $T_1$, $T_2)$, indeed, 
 $AB$ corresponds to the edge   $vi-ii$, $T_1$ rooted at $A$ ($vi$) corresponds to the uncrossed edges with endpoints on arc $v-vii$ connecting to e-graph $vi-ii$, and $T_2$  rooted at $B$ ($ii$) corresponds to 
 the uncrossed edges with endpoints on arc $i-iii$  connecting to e-graph $vi-ii$. 
Indeed,  there is a one-to-one correspondence between  $L=(AB$, $T_1$, $T_2)$ with $k$ edges and a u-c-subgraph of the u-c-tree $C$ reducing to ${T_3}^6[5]$ consisting of 
$vi-ii$,
and  all the uncrossed edges connecting to e-graph $vi-ii$ with the restriction that they have both of their endpoints on arcs $i-iii$ and $v-vii$ and that there is at least on edge on arc $ii-iii$, using 
 $k$ edges. Since we are interested in the parity of the number of  the latter with $k=\frac{n-2}{4}$, it suffices to investigate the parity of the number of
 $L=(AB$, $T_1$, $T_2)$ with $\frac{n-2}{4}$ edges in order to obtain an answer to our original question. This is the problem that we  solve in the following sections.

\begin{definition}
Let $\mathcal{A}_k=\{T \in \mathcal{A}$ $\mid T$ has $k$ edges$\}$.
Let $a_k=\mid\mathcal{A}_k\mid$.

Let $\mathcal{B}_k=\{ T \in \mathcal{B}$  $\mid T$ has $k$ edges$\}$.
Let $b_k=\mid\mathcal{B}_k\mid$. 

Let $\mathcal{L}_k=\{ L \in \mathcal{L}$  $\mid L$ has $k+1$ edges$\}$.
Let $l_k=\mid\mathcal{L}_k\mid$. 
\end{definition}

In terms of the just introduced symbols we are looking for the parity of $l_s$, where $s=\frac{n-2}{4}-1$.

\begin{lemma}
 $l_s= \sum^{s-1}_{i=0} a_{i}\cdot b_{s-i}$.
\end{lemma}
\begin{proof}
By definition, any $L\in \mathcal{L}$ is a triple $(AB$, $T_1$, $T_2)$, with a fixed edge $AB$ and $T_1\in \mathcal{A}$ rooted at $A$,  $T_2\in \mathcal{B}$ rooted at $B$. In order to obtain $l_s$ we have to sum over all possible $T_1$ and $T_2$, such that the sum of the number of edges of these two forms is $s$. Since $T_2$ has at least one edge by definition, the number of edges of $T_1$ can vary from $0$ to $s-1$. Thus,  $l_s= \sum^{s-1}_{i=0} a_{i}\cdot b_{s-i}$.\end{proof}

%In our case $s$ will equal  $\frac{n-2}{4} -1$.

\begin{definition}
A left-right tree T is a finite set of vertices such that:

a. One specially designated vertex is called the root of T and a  left-right delimiter $k\in \mathbb{N} \backslash \{0\}$ is specified,

b. The remaining vertices (excluding the root) are put into an ordered  partition $(T_1,\ldots, T_l)$  of $l\geq 0$ disjoint non-empty sets $T_1,\ldots, T_l$, each of which is a left-right tree.  The left-right delimiter $k$ specifies that
 $T_1,\ldots, T_{k-1}$ are left from the root of $T$, and the edges connecting the root of $T$ with the roots of  $T_1,\ldots, T_{k-1}$ are called left edges, while the roots of  $T_k,\ldots, T_l$ are said to be right from the root of $T$ and the edges connecting the root of $T$ and the roots of  $T_k,\ldots, T_l$   are called right edges. The edge connecting the root of T with the root of $T_l$ is called the rightmost edge of the l-r-tree.  The trees $T_1,\ldots, T_l$ are called subtrees of the root, more precisely trees $T_1,\ldots, T_{k-1}$ are left  subtrees 
while trees $T_k,\ldots, T_l$ are right subtrees of the root. 

\end{definition}

For some simple examples of l-r-trees see Figure 11.2.

There is a straightforward bijection between the set  $\mathcal{A}_k$  and the set of left-right-trees with $k$ edges. Also, 
there is a straightforward bijection between the set  $\mathcal{B}_k$  and the set of left-right-trees with $k$ edges which have at least one right edge coming out of the root.

Namely, given $T\in \mathcal{A}$ ($\mathcal{B}$) with axis $AB$, rooted at $A$ ($B$), 
let $C_1,\ldots, C_{k-1}$ be all the vertices of $T$, such that $AC_i$ ($BC_i$) is an edge, and $C_i$ is in arc $\widehat{BA}$ ($\widehat{AB}$). Let $C_k,\ldots, C_l$ be all the vertices of $T$, such that $AC_i$ ($BC_i$) is an edge, and $C_i$ is in arc $\widehat{AB}$ ($\widehat{BA}$). Arcs $AC_l$, $C_lC_{l-1}, \cdots, C_2C_1, C_1A$ ($BC_l$, $C_lC_{l-1}, \cdots, C_2C_1, C_1B$)
cover the circle and are disjoint. 
Let $T_1,\ldots, T_k$ be elements 
 of  $\mathcal{A}$, $T_i$ rooted at vertex $C_i$, having axis $AC_i$ ($BC_i$) and containing all the edges of $T$ which they can possibly contain.   

Then, the recursive definition of the bijection is:

a. Set vertex $A'$ ($B'$) to be the root of the corresponding left-right-tree $T'$ (l-r-tree for short).

b. $({T'}_1,\ldots, {T'}_l)$ are the ordered subtrees of the root of $T'$ with property that $({T'}_1,\ldots, {T'}_{k-1})$ are the left subtrees, $({T'}_k,\ldots, {T'}_l)$ are the right subtrees, and ${T'}_i$ is the l-r-tree corresponding to $T_i$. 
Finally, a single point corresponds to a single point. 

From the bijections above we deduce that the number of l-r-trees with $k$ edges is $a_k$, while the number of l-r-trees with $k$ edges such that there is at least one right edge coming out of the root is $b_k$. 
For convenience we define the following three sets:

\begin{definition}
$\mathcal{A'}_k=\{T \mid T$ is a   l-r-trees with $k$ edges$\}$.

$\mathcal{B'}_k=\{T \mid T$ is a   l-r-trees with $k$ edges  such that there is at least one right edge coming out of the root of $T\}$.

$\mathcal{C'}_k=\mathcal{A'}_k$$\setminus$$\mathcal{B'}_k$. Let $c_k=\mid\mathcal{C'}_k\mid$

\end{definition}
Now we state the results mentioned above using these symbols. We also include the relation $c_k=a_k-b_k$ which follows directly from the definition of $\mathcal{C'}_k$.

\begin{lemma}
$\mid\mathcal{A'}_k\mid=a_k$, 
$\mid\mathcal{B'}_k\mid=b_k$, and
$a_k=c_k+b_k$.

\end{lemma}

\section{Results About $a_k$, $b_k$ and $c_k$} 
 
As mentioned in the previous section, 
we are interested in the parity of the sum 
$l_s=\sum^{s-1}_{i=0} a_{i}\cdot b_{s-i}$, where $s=\frac{n-2}{4} -1$. In order to determine this parity,  
we first investigate the parities of $a_k$ and $b_k$. 
 
 In the following text we write   $\equiv$  to mean  equivalence 
$modulo$ \textit{2}.

Observe that any $T\in \mathcal{B'}_m$ consists of the followings:

$\bullet$ a right edge $B'C'$ coming out of the root $B'$ of the l-r-tree $T$; 
let  $B'C'$ be the  rightmost edge of $T$
 
$\bullet$ a $T_1 \in \mathcal{A'}_l$ rooted in $C'$

$\bullet$ a $T_2 \in \mathcal{A'}_{m-l-1}$ rooted in $B'$.

Thus,
 $b_{m}=\sum^{k-1}_{l=0} a_{l} \cdot a_{m-l-1}$. 
This equation immediately shows that if $m=2k$,  
then $b_{m} \equiv 0$, and if $m=2k+1$, then  
                             $b_{m} \equiv a_{k}^2 \equiv a_{k}$. That is:

\begin{lemma}
 $b_{2k} \equiv 0$ and 
                             $b_{2k+1} \equiv a_{k}$. 
\end{lemma}

%Note that from the definition 
%of object $C$ we have $A_k=B_k+C_k$. 

\begin{lemma} 
$ a_{2k+1}\equiv 0 $ 
\end{lemma} 

\begin{proof}                     We  prove the statement by induction. 
 
$\bullet$  Base of induction: $a_{1}=2\equiv 0$

$\bullet$ Inductive hypothesis: for all $0\leq k'<k$,   $a_{2k'+1}\equiv0$

$\bullet$ Inductive step:

Enumerate $a_{2k+1}$ summing over all $a_{2k+1, m}$, 
$m=1,2,\ldots, 2k+1$, where $a_{2k+1, m}$ is the number 
of $T\in  \mathcal{A'}_{2k+1}$ such that there are exactly
$m$ edges coming out of  the root $A'$ of $T$. 
Let  $T_1$, $T_2$, $\ldots$, $T_m$ be the subtrees of $T$, where 
 $T_i$ is the l-r-tree rooted 
at the $i^{th}$ child of $A'$ (so that the $m$ edges coming out of the root and  $T_1$, $T_2$, $\ldots$, $T_m$ have no common edges, but their union is the whole l-r-tree $T$). Let $t_i$ be the number 
of edges in $T_i$, for all $i\in \{1,2,\ldots, m\}$. Then,  
$t_1+t_2+\cdots+t_m=2k+1-m$. 
We obtain $a_{2k+1, m}$ by summing $(m+1)\cdot \prod^{m}_{i=1} a_{t_i}$ over all possible choices of $t_1,$ $\ldots,$ $t_m$, since  $a_{2k+1, m}$ is equal to the number of ways to choose which of the $m$ edges coming out of the root are right or left ($m+1$ ways for this) times the number of ways to construct the subtrees rooted at the children of the root ($\prod^{m}_{i=1} a_{t_i}$), and this all over the possible $t_1,$ $\ldots,$ $t_m$, satisfying $t_1+t_2+\cdots+t_m=2k+1-m$.

$First$ $Claim$. $a_{2k+1, 2l+1}\equiv 0$  
 
Regardless of how we fix the $t_i$ satisfying 
$t_1+t_2+\cdots+t_{2l+1}=2k+1-(2l+1)$, the sum  $(2l+2)\cdot \prod^{2l+1}_{i=1} a_{t_i}$ over all possible choices of $t_1,$ $\ldots,$ $t_m$
will be divisible by $2l+2$, 
 therefore, 
 $a_{2k+1, 2l+1}$ is  divisible by $2$.

$Second$ $Claim$. $a_{2k+1, 2l}\equiv 0$

There are $2l$ edges coming out of the root $A'$, thus, 
$t_1+t_2+\cdots+t_{2l}=2k+1-2l\equiv 1$.
Therefore, in every case when there are $2l$ edges coming out of the root   at least one $t_j$ must be  odd, thus there exists a $j$ such that  $a_{t_j} \equiv 0$ by the inductive hypothesis 
(since $2l\geq2$ it follows that  $t_j<2k+1$, so we can use the inductive hypothesis).
Since $a_{2k+1, 2l}$ is the sum of 
$(2l+1)\cdot \prod^{2l}_{i=1} a_{t_i}$ 
 over all possible 
 $t_1$, $t_2,$ $\ldots,$ $t_{2l}$,  
%i.e.  $A^{1},A^{2},\ldots,A^{2l}$, 
 then 
 this sum is $\equiv 0$, having that in each product there is some  
$a_{t_j} \equiv 0$.

Since $a_{2k+1}=\sum^{2k+1}_{m=1} a_{2k+1,m}$, and $a_{2k+1,m}\equiv 0$ 
for all $m$, then $a_{2k+1}\equiv 0$, and the induction is finished. 
\end{proof} 
 
\begin{lemma} 
$a_{2k}\equiv c_{k}$ 
\end{lemma} 

\begin{proof} 
Let $F_1$ be a  l-r-tree with $2k$ edges.
Let $R$ be the root of $F_1$, let $d$ be its delimiter, and let  $(T_1,\ldots, T_l)$ be the l-r-subtrees rooted at the children of the root ($T_i$ is the l-r-tree having the $i^{th}$ child for its root).  We define operation $f$, ``flip,'' as follows. Let l-r-tree $F_2=f(F_1)$  have root $R$, delimiter $n-d+2$ and the ordered partition of the remaining vertices is  $({T'}_l,\ldots, {T'}_1)$, where ${T'}_i=f(T_i)$. If $F_1$ contains just one vertex, and no edges, then $f(F_1)=F_1$.  Intuitively, flipped $F_1$ is nothing but the l-r-tree flipped over a vertical axis going through the root, Figure 11.1. It is clear that if $f(F_1)=F_2$ then $f(F_2)=F_1$. 
\begin{figure}
\begin{center}
\epsfbox{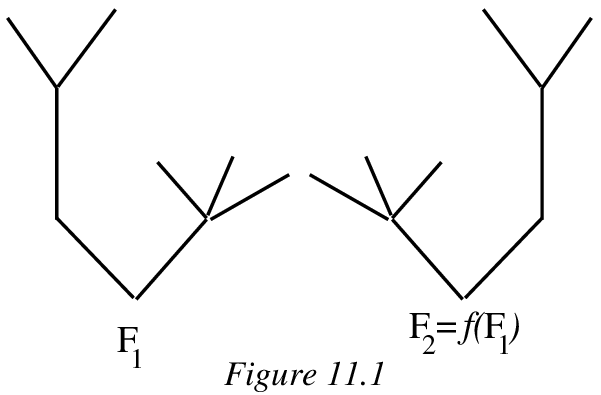}
\end{center}
\end{figure}
We consider two l-r-trees identical if the   
edge connecting the root and the $i^{th}$ child is left in both l-r-trees or right in both l-r-trees  and if the l-r-subtrees rooted each child are identical.

Pairs of l-r-trees ($F_1$, $F_2$), where $f(F_1)=F_2$, and  $F_1$ and $F_2$  are not identical
contribute an even number to the number of l-r-trees with $2k$ edges, and so, 
the parity of $a_{2k}$ is the parity of the number of l-r-trees $F$  
such that $f(F)=F$.
In order for $f(F)$ to be identical to $F$, 
it must be that it is ``symmetric,'' that is, the number of left edges coming out of the root  is equal to the number of rigft edges coming out of the root, all subtrees rooted at children of the root are symmetric, and finally, a point is symmetric, Figure 11.2. The subgraph of a symmetric $F$ which contains all the left edges coming out of the root and all the edges of the subtrees rooted at these left children of the root uniquely determine a symmetric $F$, furthermore, there is an obvious  one-to-one correspondence between these subgraphs and symmetric l-r-trees. Also, there is an obvious correspondence between the subgraph described and elements of $\mathcal{C'}_k$. Thus, 
the parity of $a_{2k}$ is the parity of $\mid\mathcal{C'}_k\mid$. 
Thus, $a_{2k}\equiv c_{k}$. 
\begin{figure} 
\begin{center}
\epsfbox{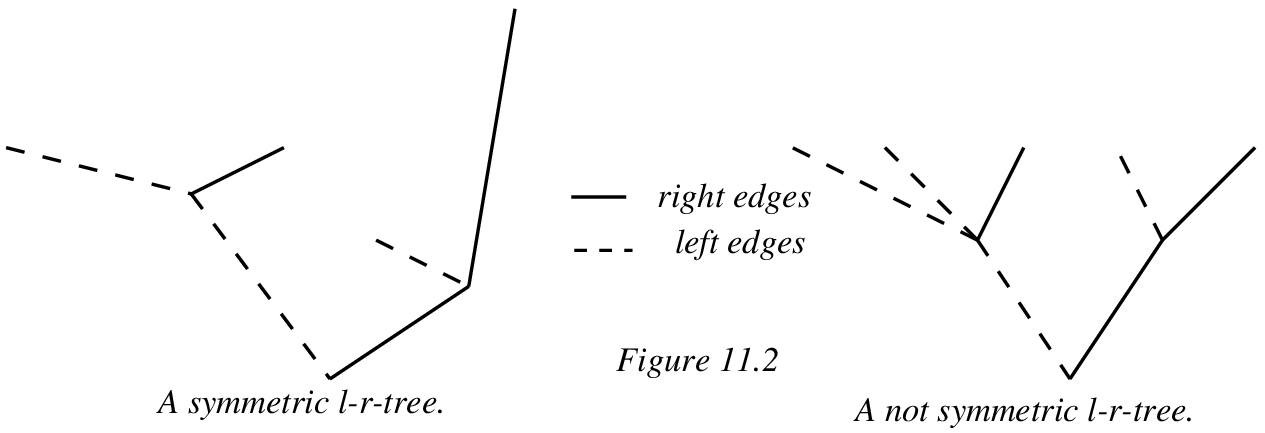} 
 
 \end{center}
\end{figure}
\end{proof} 
 
The results obtained so far:  
$i)$ $b_{2k} \equiv 0$; 
$ii)$ $b_{2k+1} \equiv a_{k}$;  
$iii)$ $a_{k}=b_{k}+c_{k}$; 
$iv)$  $a_{2k+1}\equiv0$; 
$v)$  $a_{2k}\equiv c_{k}$; 
$vi)$ From $ii)$  $b_{4k+1}\equiv a_{2k}$ and $b_{4k+3}\equiv a_{2k+1}$; 
$vii)$ Using $vi)$ and $iv)$ $b_{4k+3}\equiv 0$;
$viii)$ From $iii)$ and $i)$ $a_{2k}=b_{2k}+c_{2k}\equiv c_{2k}$, 
 that is $a_{2k}\equiv c_{2k}$; 
$ix)$ From $iii)$ and $vi)$ $a_{4k+1}=b_{4k+1}+c_{4k+1}$ and so 
 $a_{4k+1}\equiv b_{4k+1}+c_{4k+1} \equiv a_{2k}+c_{4k+1}$. Using 
 $iv)$ we get $a_{2k}\equiv c_{4k+1}$;  
$x)$ From $iii)$ and $vii)$ $a_{4k+3}\equiv b_{4k+3}+c_{4k+3}\equiv c_{4k+3}$. 
Using $iv)$ $c_{4k+3}\equiv 0$.
In all cases $k\geq 0$. 
 
Summarizing this:

\textbf{$a)$ $a_{2k+1}\equiv 0 $; $b)$ $a_{2k} \equiv c_{4k+1} \equiv c_{2k} \equiv c_{k}$; $c)$ $c_{4k+3}\equiv 0$; $d)$ $b_{2k} \equiv 0$; $e)$ $b_{4k+3}\equiv 0$; $f)$ $b_{4k+1}\equiv a_{2k}$}
 In all cases $k\geq 0$. 
 
We now directly attack the problem of determining the parity of  $l_s=\sum^{s-1}_{i=0} a_{i}\cdot b_{s-i}$, where 
$s= \frac{n-2}{4}-1$, $s\geq 1$.

\textit{\textbf{Case 1.}} If $s=2r$, for some $r\geq 1$, either $i$ is odd, or $s-i$ is even, therefore 
 $\sum^{s-1}_{i=0} a_{i}\cdot b_{s-i} \equiv 0$ by 
  $a)$  and $d)$. 
 Therefore, $l_s$ is even for $\frac{n-2}{4}-1=s=2r$, where $r\geq 1$.

\textbf{\textit{Result.}} For $n=8r+6$, $r\geq 1$, $f(n)$ is divisible by $n$.

\textit{\textbf{Case 2.}} If  $s=2r+1$, for some $r\geq 0$, $l_s$  is equivalent 
   with the sum over those $i$ for which $i=2k$ and $s-i=4l+1$, since
  $a_{2k+1}\equiv 0$, $b_{2k} \equiv 0$, and  $b_{4k+3} \equiv 0$. We analyze this case by splitting it into more smaller 
cases.  
 
\textit{\textbf{Case 2.1.}}: $s=4v+1$, $2k+4l+1=4v+1$, $v,k,l\geq 0$

Then $k=2(v-l)$. Since $l$ can be $0,1,2,\ldots,v$; 
$k$ is then $2v, 2v-2,\ldots,0$. 
 
Using $b)$ and $f)$ we have  $a_{2k} \equiv c_{k}$ and $b_{4k+1} \equiv c_{k}$,
and so 
  $\sum^{s-1}_{i=0} a_{i}\cdot b_{s-i}\equiv 
  \sum^v_{i=0} c_{2(v-i)}\cdot c_{i}\equiv$ 
 $ \sum^v_{i=0} c_{i}\cdot c_{2(v-i)}$. 
Using  $c_{2k}\equiv c_{k}$, 
$\sum^{s-1}_{i=0} a_{i}\cdot b_{s-i}\equiv \sum^v_{i=0} c_{i}\cdot c_{v-i}$. 
The last sum is symmetric, thus in case $v=2v_1+1$, for some $v_1\geq 0$,  this sum is even, that is $l_s$ is even
for 
 $\frac{n-2}{4}-1=s=4v+1=4(2v_1+1)+1=8v_1+5$. 

\textbf{\textit{Result.}} For $n=32v_1+26$, $v_1\geq 0$, $f(n)$ is divisible by $n$.

In the  case when $v=2v_1$, $v_1\geq 0$ we have:

\textbf{\textit{Result.}}
For $n=32v_1+10$, $v_1\geq 0$ :
$l_s\equiv c_{v_1}$, $s=\frac{n-2}{4}-1$.

\textit{\textbf{Case 2.2.}}:  $s=4v+3$, $2k+4l+1=4v+3$, $v,k,l\geq 0$

Then $k=2(v-l)+1$. Since $l$ can be $0,1,2,\ldots,v$; 
$k$ is then $2v+1, 2v-1,\ldots,1$. 
Using $b)$ and $f)$ we have  $a_{2k} \equiv c_{k}$ and $b_{4k+1} \equiv c_{k}$, 
and so 
  $\sum^{s-1}_{i=0} a_{i}\cdot b_{s-i}\equiv \sum^v_{i=0} c_{i}\cdot c_{2(v-i)+1}$. 
  
\textit{\textbf{Case 2.2.1.}}:  $v=2v_1$,  $v_1\geq 0$. 
 
Using $c_{4k+3}\equiv 0$ we get 
$\sum^{s-1}_{i=0} a_{i}\cdot b_{s-i}\equiv\sum^{v_1}_{i=0} c_{2i}\cdot c_{2(v-2i)+1}$. 
Using $c_{2k}\equiv c_{k}$ we get 
$\sum^{s-1}_{i=0} a_{i}\cdot b_{s-i}\equiv\sum^{v_1}_{i=0} c_{i}\cdot c_{2(v-2i)+1}$. 
Since $2(v-2i)+1=2(2v_1-2i)+1=4(v_1-i)+1$, and $c_{4k+1}\equiv c_{k}$, we get  
$\sum^{s-1}_{i=0} a_{i}\cdot b_{s-i}\equiv\sum^{v_1}_{i=0} c_{i}\cdot c_{v_1-i}$. 
Since $\sum^{v_1}_{i=0} c_{i}\cdot c_{v_1-i}$ is symmetric, we have the following 
result: 
for $v_1=2v_2+1$,  $v_2\geq 0$, the sum is even, and so  $l_s$ is even for 
$\frac{n-2}{4}-1=s=4v+3=4\cdot2v_1+3=8(2v_2+1)+3=16v_2+11$.

\textbf{\textit{Result.}} For
$n=64v_2+50$,  $v_2\geq 0$, 
 $f(n)$  is divisible by $n$.

In the case when $v_1=2v_2$ we have:

\textbf{\textit{Result.}} For $n=64v_2+18$, $v_2\geq 0$:   $l_s\equiv \sum^{v_1}_{i=0} c_{i}\cdot c_{v_1-i}\equiv c_{v_2}$.

\textit{\textbf{Case 2.2.2.}}: $v=2v_1+1$,  $v_1\geq 0$. 
 
Using $c_{4k+3}\equiv 0$ from 
$\sum^{s-1}_{i=0} a_{i}\cdot b_{s-i}\equiv 
\sum^v_{i=0} c_{i}\cdot c_{2(v-i)+1}$ 
we get   
$\sum^{s-1}_{i=0} a_{i}\cdot b_{s-i}\equiv 
\sum^{v_1}_{i=0} c_{2i+1}\cdot c_{2(v-2i-1)+1}\equiv 
\sum^{v_1}_{i=0} c_{2i+1}\cdot c_{2(2v_1-2i)+1} \equiv
\sum^{v_1}_{i=0} c_{2i+1}\cdot c_{4(v_1-i)+1}$ 
and using $c_{4k+1}\equiv c_{k}$ we have 
$\sum^{s-1}_{i=0} a_{i}\cdot b_{s-i}\equiv 
\sum^{v_1}_{i=0} c_{2i+1}\cdot c_{4(v_1-i)+1}\equiv 
\sum^{v_1}_{i=0} c_{2i+1}\cdot c_{v_1-i}$. 
Since, $\sum^{v_1}_{i=0} c_{2i+1}\cdot c_{v_1-i}= 
\sum^{v_1}_{i=0} c_{v_1-i}\cdot c_{2i+1}= 
\sum^{v_1}_{j=0} c_{j}\cdot c_{2(v_1-j)+1}$, we get that 
$\sum^{s-1}_{i=0} a_{i}\cdot b_{s-i}\equiv 
\sum^{v_1}_{i=0} c_{i}\cdot c_{2(v_1-i)+1}$. 
Note that we started from  
$\sum^{s-1}_{i=0} a_{i}\cdot b_{s-i}\equiv 
\sum^v_{i=0} c_{i}\cdot c_{2(v-i)+1}$, and 
sums $\sum^v_{i=0} c_{i}\cdot c_{2(v-i)+1}$ and 
$\sum^{v_1}_{i=0} c_{i}\cdot c_{2(v_1-i)+1}$ are of the same form, 
with the difference that the second one is of length $v_1+1$, while the first 
 is of length $v+1=2v_1+2$, that is  twice the second. 

Let $S(v)=\sum^v_{i=0} c_{i}\cdot c_{2(v-i)+1}$.
 
 Observe, that we got 
$l_s\equiv S(v) \equiv S(v_1)$
 and we still 
 do not know what this sum is equivalent to modulo 2. Observe, that 
since the last two sums are of the same form, we can now begin the 
process of \textit{Case 2.2.}, that is:  
 
\textit{Case 2.2.1.}:  $v_1=2v_2$, if $v_2=2v_3+1$ the sum is even, 
 but if  $v_2=2v_3$ the sum is $\equiv c_{v_3}$ 
 
\textit{Case 2.2.2.}: $v_1=2v_2+1$, then 
$\sum^{s-1}_{i=0} a_{i}\cdot b_{s-i}\equiv 
\sum^v_{i=0} c_{i}\cdot c_{2(v-i)+1}\equiv 
 \sum^{v_1}_{i=0} c_{i}\cdot c_{2(v_1-i)+1}\equiv 
  \sum^{v_2}_{i=0} c_{i}\cdot c_{2(v_2-i)+1}$ 
and so forth. Depending on $v_2$ we either continue the process, or 
get the final result in case $v_2$ is even. 
 
Since in case $v$, $v_1$, $v_2$, etc. are odd, the sums $S(v)$, $S(v_1)$, $S(v_2)$, etc.  we are considering 
always get twice shorter (from $v+1$ to $v_1+1$ to $v_2+1$, etc.) and so 
this process is finite (we cannot half an integer infinitely many times and 
still get an integer), i.e. at one point either the $v_k$ we are going 
to consider must be even and greater than $0$, and then we apply 
the \textit{Case 2.2.1.} and get 
either that the sum is even, or that the sum $\equiv c_{v_k}$, or  $v_k$ 
  gets equal to \textit{1}, and in this case the sum is odd.

\textbf{\textit{Result.}}
For $n=32v_1+34$, $v_1\geq 0$, we get a sum $S(v_1)$ for which we have to 
decide depending on $v_1$ which subcase of $\textit{Case 2.2}$  applies.

\section{The Behavior of  $c_{k}$}

In order to clarify the situation in the case where $n=32v_1+10$, $v_1\geq 0$, and $n=64v_2+18$, 
$v_2\geq 0$, we have to determine how $c_{k}$ behaves. We use the facts that  $c_{4k+1}\equiv c_{2k}\equiv c_{k}$  and 
$c_{4k+3}\equiv 0$. 
 
Let an integer $v>0$ of the form $4k+3$, $4k+1$ or $2k$ be given.  

Let $V=v$. Therefore $c_{V}\equiv c_{v}$.  
 
\textit{\textbf{Step I}}: 
 
if $v=4k+3$, we know that  $c_{V}\equiv c_{v}\equiv 0$: proceed to \textit{Step II}   
 
if $v=4k+1$, $v>1$, we know that $c_{V}\equiv c_{v} \equiv c_{4k+1}\equiv c_{k}$ and $k<4k+1$:  redefine $v:=k$ and   proceed to \textit{Step II}

if  $v=2k$, $v>1$, we know that  $c_{V}\equiv c_{v}\equiv c_{2k}\equiv c_{k}$ and $k<2k$: redefine $v:=k$ and   proceed to \textit{Step II}

if $v=1$, we know that  $c_{V}\equiv c_{v}\equiv 1$: proceed to \textit{Step II}

\textit{\textbf{Step II}}: 

continue \textit{Step I} until $v=4k+3$ for some $k$, that is $c_{V}\equiv 0$, or 
$v=1$, that is  $c_{V}\equiv 1$.

It is clear that by executing the process described above, $v$ must 
become either $1$, or of form $4k+3$ in a  finite number of steps. 
 
Since only those $c_{v}$, for which $v$ gets to 1 are odd, we aim for 
determining the form of these $v$.

Define: 
 
$\bullet$ unary operation $m$ on an integer $k$: $m(k)=2\cdot k$ 
 
$\bullet$ unary operation $M$ on an integer $k$: $M(k)=4\cdot k+1$ 
 
It follows form the process described above that all integers $v$ for which $c_{v}\equiv 1$ can be created by any 
 number and any order of operations $m$ and $M$ on $1$, e.g.

$ v =M\circ M\circ m\circ m\circ m\circ M\circ m\circ m (1)=$

 $   M\circ M\circ m\circ m\circ m\circ M\circ m (2)=$

$    M\circ M\circ m\circ m\circ m\circ M (4)=         $

$    M\circ M\circ m\circ m\circ m (17)=$

$    M\circ M\circ m\circ m (34)=         $

$    M\circ M\circ m (68)=$

$    M\circ M (136)=        $

 $   M (545)= 2181$ 
 
$\Rightarrow$  $c_{2181}\equiv 1$, thus, for example, for 
$n=32\cdot 2181+10=69802$ 
the number of genus one c-trees on $69802$ points is not divisible 
 by $69802$, only by $69802\div 2=34901$. 
 
It is easy to see, that all $v$ obtained by finitely many operations 
$M$ and $m$ on 1 are of the following form: 
 
$1\cdot 2^{l_0}+4\cdot 2^{l_1}+4^2\cdot 2^{l_2}+4^3\cdot 2^{l_3}+ 
\cdots+4^i\cdot 2^{l_i}+\cdots+4^f\cdot 2^{l_f}$, where 
 
$0\leq l_0\leq l_1\leq l_2\leq l_3\leq \cdots \leq l_i\leq \cdots\leq l_f$, $l_i \in \mathbb{N}$, $i \in \{0, 1, 2,\ldots, f\}$.

If $v$ satisfies the previous or if $v=0$ ($c_0 \equiv 1$), 
we will say $v$ is a \textbf{\textit{negligent number}} (or $nn$ for short). 
 
Therefore, we have proven that   $(c_{v}\equiv 1)$ $\Leftrightarrow$   $(v$ is $nn)$.  
For all other $v$, $c_{v}\equiv 0$. 
From this we have that:

\textbf{\textit{Results.}} For $n=32v_1+10$, $v_1\geq 0$, where $v_1$ is a $nn$ 
$l_s$ is odd, that is $f(n)$ is not divisible by $n$ but is divisible by $\frac{n}{2}$.
 For $n=32v_1+10$, $v_1\geq 0$, where $v_1$ is not a  $nn$ 
$l_s$ is even, that is $f(n)$  is divisible by $n$. 
 
For $n=64v_2+18$, $v_2\geq 0$, where $v_2$ is  a $nn$
$l_s$ is odd, that is $f(n)$ is not divisible by $n$ but is divisible by $\frac{n}{2}$. 
 For $n=64v_2+18$, $v_2\geq 0$, where $v_2$ is not a $nn$
$l_s$ is even, that is $f(n)$ is divisible by $n$.

\section{What If $n=32k+34$?} 
 
In Section 11 we have seen that for $n=32v_1+34$, $v_1\geq 0$, we get a sum $S(v_1)$ for which we have to 
decide depending on $v_1$ which subcase of  \textit{Case 2.2}  applies in order to determine the parity of $l_s$. We decide this now. 
 
We have that 
$\frac{n-2}{4}-1=s=4v+3$ and  $v=2v_1+1$, $v_1\geq 0$. 
There are two possibilities: 
 
%1. Possibility 
 
%$v=v_0=2v_1+1$, $v_1=2v_2+1$, $v_2=2v_3+1$,$\ldots$, $v_{k-1}=2v_{k}+1$, $v_{k}=1$, 
%$k\geq 0$. 
%It is easy to see, that $v=1+2+2^2+\cdots+2^{k}=2^{k+1}-1$, thus 
%$s=4\cdot 2^{k+1}-1$, i.e. $n=16\cdot 2^{k+1}+2$, $k\geq 0$. 
%In this case the sum is odd. 
 
$First$ $Possibility$: 
 
$v=2v_1+1$, $v_1=2v_2+1$, $v_2=2v_3+1$,$\ldots$, $v_{k-1}=2v_{k}+1$, 
$v_{k}=2v_{k+1}$, $v_{k+1}=2v_{k+2}+1$, $k\geq 1$. 
It is easy to see, that $v=1+2+2^2+\cdots+2^{k-2}+2^{k-1}+2^{k} \cdot v_k 
=2^{k}-1+2^{k} \cdot v_k$, thus 
$s=16\cdot 2^{k}\cdot v_{k+2}+3\cdot4\cdot2^{k}-1$, 
that is  $n=64\cdot 2^{k}\cdot v_{k+2}+3\cdot16\cdot2^{k}+2$, $k\geq 1$. 
In this case $l_s$ is even.

$Second$ $Possibility$: 
 
$v=2v_1+1$, $v_1=2v_2+1$, $v_2=2v_3+1$,$\ldots$, $v_{k-1}=2v_{k}+1$, 
$v_{k}=2v_{k+1}$, $v_{k+1}=2v_{k+2}$,  $k\geq 1$, $v_{k+2}>0$.  
It is easy to see, that $v=1+2+2^2+\cdots+2^{k-2}+2^{k-1}+2^{k} \cdot v_k 
=2^{k}-1+2^{k} \cdot v_k$, thus 
$s=16\cdot 2^{k}\cdot v_{k+2}+4\cdot2^{k}-1$, 
that is $n=64\cdot 2^{k}\cdot v_{k+2}+16\cdot2^{k}+2$, $v_{k+2}>0$. 
In this case  $l_s \equiv c_{v_{k+2}}$, and  including the case 
when $v_{k+2}=0$  we get that  
if $v_{k+2}$ is a $nn$ 
 then $l_s$ is odd, otherwise $l_s$ is even. 
 
Now that we have obtained the parity of $l_s$ in every case, we can state the main theorem of our work, namely when $f(n)$ is divisible by $n$ and when it is divisible only by $\frac{n}{2}$.
 
\section{The Theorem}

\textit{The number of genus one l-c-trees on $n$ points  is divisible by $n$ or  $\frac{n}{2}$ for all integers $n>3$.}

\textit{The number of genus one l-c-trees on $n$ points ($n>3$) is not divisible by $n$, but  is  divisible by $\frac{n}{2}$  if and only if $n$ is of form:}

$32v_1+10$, $v_1\geq 0$, and $v_1 $ is a negligent number; $64v_2+18$, $v_2\geq 0$, and $v_2$ is  a negligent number; 
     $64\cdot 2^{k}\cdot v_{k+2}+16\cdot2^{k}+2$,$k\geq 1$, and $v_{k+2}$ 
   is a negligent number.

 \textit{For all other $n$ the number of genus one l-c-trees on $n$ points ($n>3$) is divisible by $n$.}
 
%The sum will be odd for $n$ of form: 

\begin{proof} 
 
The proof is given in sections 1 to 13. 
\end{proof}

\section*{Acknowledgments} 
 
I am  grateful to Professor Richard P. Stanley for his encouragement and advice  throughout this research.  
 I  would also like to thank Thomas Lam who read and commented on drafts of this paper. 
This research was supported by the UROP Office and the Lord Foundation.

\begin{singlespace}

\end{singlespace} 

\end{singlespace}

\end{document}